\documentclass[11pt]{article}
\usepackage[utf8]{inputenc}
\usepackage[english]{babel}
\usepackage{hyperref}
\usepackage{amsmath,amsfonts,amscd,amssymb}
\usepackage[all]{xy}
\usepackage{makeidx}
\usepackage{color}
\usepackage{graphicx}
\usepackage[active]{srcltx} 
\newtheorem{thm}{Theorem}[section]
\newtheorem{cor}[thm]{Corollary}
\newtheorem{lem}[thm]{Lemma}
\newtheorem{prop}[thm]{Proposition}
\newtheorem{deft}[thm]{Definition}
\newtheorem{rek}[thm]{Remark}
\newtheorem{conj}{Conjecture}
\let\noi=\noindent

\let\sur=\overline
\let\sous=\underline

\newcommand{\dsp}{\displaystyle}
\def\N{\mathbb{N}} 
\def\Z{\mathbb{Z}} 
\def\R{\mathbb{R}} 
\def\C{\mathbb{C}} 
\def\Q{\mathbb{Q}} 
\def\F{\mathbb{F}}
\def\sper{\mbox{Sper}}

\def\gr{\mbox{gr}}
\def\init{\mbox{in}}
\def\notin{\mbox{$\in$ \hspace{-.8em}/}} 
\def\Ra{\Rightarrow} 
\def\Lra{\Leftrightarrow} 
\def\sgn{\mbox{sgn}} 
\makeindex
\title{ON THE PIERCE--BIRKHOFF CONJECTURE}
\author{F. Lucas\\LUNAM, Universit{\'e} d'Angers\\
LAREMA CNRS UMR 6093\\
2, bd Lavoisier\\
49045 Angers C{\'e}dex, France \and
D. Schaub\\LUNAM, Universit{\'e} d'Angers\\
LAREMA CNRS UMR 6093\\
2, bd Lavoisier\\
49045 Angers C{\'e}dex, France \and
M. Spivakovsky\\
Institut de Math\'ematiques de Toulouse/CNRS UMR 5219\\
Universit\'e Paul Sabatier\\
118, route de Narbonne\\
31062 Toulouse cedex 9, France.}

\date{}

\begin{document}
\maketitle


\abstract{This paper represents a step in our program towards the proof of the
Pierce--Birkhoff conjecture. In the nineteen eighties J. Madden proved that the
Pierce-Birkhoff conjecture for a ring $A$ is equivalent to a statement about an
arbitrary pair of points $\alpha,\beta\in\sper\ A$ and their separating ideal
$<\alpha,\beta>$; we refer to this statement as the \textbf{local
Pierce-Birkhoff conjecture} at $\alpha,\beta$. In \cite{LMSS} we introduced a
slightly stronger conjecture, also stated for a pair of points
$\alpha,\beta\in\sper\ A$ and the separating ideal $<\alpha,\beta>$, called the
\textbf{Connectedness conjecture}. In this paper, for each pair
$(\alpha,\beta)$ with $ht(<\alpha,\beta>)=\dim A$, we define a natural
number, called complexity of $(\alpha,\beta)$. Complexity 0 corresponds to the
case when one of the points
$\alpha,\beta$ is monomial; this case was settled in all dimensions in
\cite{LMSS}. In the present paper we introduce a new conjecture, called the
\textbf{Strong Connectedness conjecture}, and prove that the strong
connectedness conjecture in dimension $n-1$ implies the connectedness conjecture
in dimension $n$ in the case when $ht(<\alpha,\beta>) \le n-1$. We prove the
Strong Connectedness conjecture in dimension 2, which gives the Connectedness
and the Pierce--Birkhoff conjectures in any dimension in the case when
$ht(<\alpha,\beta>) \le2$. Finally, we prove the Connectedness (and hence also
the Pierce--Birkhoff) conjecture in the case when $\dim\
A=ht(<\alpha,\beta>)=3$, the pair $(\alpha,\beta)$ is of complexity 1 and $A$ is
excellent with residue field $\R$.}

\section{Introduction}

All the rings in this paper will be commutative with 1. Let $R$ be a
real closed field. Let $B=R[x_1,\dots,x_n]$. If $A$ is a ring and
$\mathfrak p$ a prime ideal of $A$,
$\kappa(\mathfrak p)$ will denote the residue field of $\mathfrak p$.

The Pierce--Birkhoff conjecture asserts that any
piecewise-polynomial function $f:R^n\rightarrow R$ can be expressed as
a maximum of minima of a finite family of polynomials in $n$
variables.
We start by giving the precise statement of the conjecture as it was
first stated by M. Henriksen and J. Isbell in the early nineteen
sixties.
\begin{deft}\label{pw} A function $f:R^n\to R$ is said to be \textbf{piecewise
  polynomial} if $R^n$ can be covered by a finite collection of closed
  semi-algebraic sets $P_i$ such that for each $i$ there exists a
  polynomial $f_i\in B$ satisfying
  $\left.f\right|_{P_i}=\left.f_i\right|_{P_i}$.
\end{deft}
Clearly, any piecewise polynomial function is continuous. Piecewise
polynomial functions form a ring, containing $B$, which is denoted by
$PW(B)$.\medskip

On the other hand, one can consider the (lattice-ordered) ring of all
the functions obtained from $B$ by iterating the operations of $\sup$
and $\inf$. Since applying the operations of sup and inf to
polynomials produces functions which are piecewise polynomial, this
ring is contained in $PW(B)$ (the latter ring is closed under $\sup$
and $\inf$). It is natural to ask whether the two rings coincide. The
precise statement of the conjecture is:
\begin{conj}\textnormal{\textbf{(Pierce-Birkhoff)}}\label{PB} If
  $f:R^n\to R$ is in $PW(B)$, then there exists a finite family of
  polynomials $g_{ij}\in B$ such that
  $f=\sup\limits_i\inf\limits_j(g_{ij})$ (in other words, for all
  $x\in R^n$, $f(x)=\sup\limits_i\inf\limits_j(g_{ij}(x))$).
\end{conj}

This paper is a step in a program for proving the
Pierce--Birkhoff conjecture. The starting point of this
program is the abstract formulation of the conjecture in terms of the real
spectrum of $B$ and separating ideals  proposed by J. Madden in
1989 \cite{Mad1}.

For more information about the real spectrum, see \cite{BCR}; there is
also a brief introduction to the real spectrum and its relevance to
the Pierce--Birkhoff conjecture in the Introduction to \cite{LMSS}.
\medskip

\noi\textbf{Terminology}: If $A$ is an integral domain, the phrase
``valuation of $A$'' will mean ``a valuation of the field of fractions
 of $A$, non-negative on $A$''. Also, we will sometimes commit the
 following abuse of notation. Given a ring $A$, a prime ideal
 $\mathfrak p\subset A$, a valuation $\nu$ of $\frac A{\mathfrak p}$
 and an element $x\in A$, we will write $\nu(x)$ instead of
 $\nu(x\mod\mathfrak p)$, with the usual convention that
 $\nu(0)=\infty$, which is taken to be greater than any element of the
 value group.
\medskip

\noi \textbf{Recall some notation} : For a point $\alpha \in \sper\ A$ we
denote by $\mathfrak{p}_\alpha$ the support of $\alpha$, by $A[\alpha] =
\frac{A}{\mathfrak{p}_\alpha}$ and by $A(\alpha)$ the field of fractions of
$A[\alpha]$. We also let $\nu_\alpha$ denote the valuation associated to
$\alpha$, $\Gamma_\alpha$ the value group, $R_{\nu_\alpha}$ the valuation
ring, $k_\alpha$ its residue field and $\gr_\alpha(A)$ the graded ring
associated to the valuation $\nu_\alpha$. For $f \in A$ with
$\gamma=\nu_\alpha(f)$, let $in_\alpha f$ denote the natural image of $f$ in
$\frac{P_\gamma}{P_{\gamma_+}}$. Finally, if $k$ is any field, we denote by
$\sur{k}$ its real closure. 

\begin{deft}
Let
$$
f :\sper\ A \to \coprod_{\alpha \in \sper\ A} A(\alpha)
$$
be a map such that, for each $\alpha \in \sper\ A$, $f(\alpha) \in A(\alpha)$. We
say that $f$ is piecewise polynomial (denoted $f \in PW(A)$) if there exits a
covering of $\sper\ A$ by a finite family $(S_i)_{i \in I}$ of
constructible sets, closed in the spectral topology and a family
$(f_i)_{i \in I}$, $f_i \in A$ such that, for each $\alpha \in S_i$,
$f(\alpha)=f_i(\alpha)$.

We call $f_i$ a local representative of $f$ at $\alpha$ and denote it
by $f_\alpha$ ($f_\alpha$ is not, in general, uniquely determined by
$f$ and $\alpha$; this notation means that one such local
representative has been chosen once and for all).
\end{deft}

\begin{deft}
A ring $A$ is a Pierce-Birkhoff ring if, for each $f \in PW(A)$, there
exists a finite collection $\{f_{ij}\}_{i,j}  \subset A$ such that $f=\sup_i
\inf_j
f_{ij}$.
\end{deft}

The generalized Pierce-Birkhoff Conjecture says:

\begin{conj} (\textbf{Pierce-Birkhoff Conjecture for regular rings})
 Let $A$ be a regular ring. Then $A$ is a Pierce-Birkhoff ring.
\end{conj}

Madden reduced the Pierce--Birkhoff
conjecture to a purely local statement about separating ideals and the
real spectrum. Namely, he introduced 

\begin{deft} Let $A$ be a ring. For $\alpha,\beta\in\mbox{Sper}\ A$,
  the \textbf{separating ideal} of $\alpha$ and $\beta$, denoted by
  $<\alpha,\beta>$, is the ideal of $A$ generated by all the elements
  $f\in A$ which change sign between $\alpha$ and $\beta$, that is,
  all the $f$ such that $f(\alpha)\geq0$ and $f(\beta)\leq0$.
\end{deft}

\begin{deft} \label{def:localPB}
A ring $A$ is locally Pierce-Birkhoff at $\alpha, \beta$ if the
following condition holds : let $f$ be a piecewise  polynomial
  function, let $f_\alpha \in A$ be a local representative of $f$ at
  $\alpha$ and $f_\beta \in A$ a local representative of $f$ at
  $\beta$. Then $f_\alpha-f_\beta \in <\alpha,\beta>$.
\end{deft}

\begin{thm} (Madden)
A ring $A$ is Pierce-Birkhoff if and only if it is locally
  Pierce-Birkhoff for all $\alpha, \beta \in \sper\ A$.
\end{thm}


\begin{rek} \label{rek:spezial}
Assume that $\beta$ is a specialization of $\alpha$. Then

(1) $<\alpha,\beta> = \mathfrak{p}_\beta$.

(2) $f_\alpha -f_\beta \in \mathfrak{p}_\beta$. Indeed, we may assume
that $f_\alpha \neq f_\beta$, otherwise there is nothing to
prove. Since $\beta \in \sur{\{\alpha\}}$, $f_\alpha$ is also a local
representative of $f$ at $\beta$. Hence $f_\alpha(\beta)
-f_\beta(\beta)=0$, so $f_\alpha-f_\beta \in \mathfrak{p}_\beta$.

Therefore, to prove that a ring $A$ is Pierce-Birkhoff, it is
sufficient to verify the Definition \ref{def:localPB} for all
$\alpha,\beta$ such that neither of $\alpha,\beta$ is a specialization
of the other.
\end{rek}

In \cite{LMSS}, we introduced
\begin{conj}\textnormal{\textbf{(the Connectedness
conjecture)}}\label{conn} Let $A$ be a regular ring. Let
$\alpha,\beta\in\sper\ A$ and let
  $g_1,\dots,g_s$ be a finite collection of elements of
  $A\setminus<\alpha,\beta>$. Then there exists a connected set
  $C\subset\sper\ A$ such that $\alpha,\beta\in C$ and
  $C\cap\{g_i=0\}=\emptyset$ for $i\in\{1,\dots,s\}$ (in other words,
  $\alpha$ and $\beta$ belong to the same connected component of the
  set $\sper\ A\setminus \{g_1\dots g_s=0\}$).
\end{conj}

In the paper \cite{LMSS}, we stated the Connectedness
conjecture (in the special case when $A$ is a polynomial ring) and proved that
it implies the Pierce--Birkhoff conjecture. The same proof shows that the
Connectedness Conjecture implies the Pierce-Birkhoff Conjecture for an
arbitrary ring.

\begin{deft}\label{cd} A subset $C$ of $\sper\ A$ is said to be
  \textbf{definably connected} if it is not a union of two non-empty disjoint
  constructible subsets, relatively closed for the spectral topology.
\end{deft}

\begin{deft} \textbf{Definable Connectedness Property}
 Let $A$ be a ring. Let $\alpha,\beta \in \sper\ A $.
We say that $A$ has the Definable Connectedness Property at $\alpha,\beta$ if,
for any finite collection $g_1,\dots,g_s$ of elements of $A \setminus
<\alpha,\beta>$, there exists a definably connected set $C\subset\sper\ A$
such
that $\alpha,\beta\in C$ and $C\cap\{g_i=0\}=\emptyset$ for
$i\in\{1,\dots,s\}$
(in other words, $\alpha$ and $\beta$ belong to the same definably connected
component of the set $\sper\ A\setminus \{g_1\dots g_s=0\}$).
\end{deft}

\begin{conj} \label{DCC} \textnormal{\textbf{(Definable Connectedness
Conjecture)}} Let $A$ be a regular ring. Then $A$ satisfies the Definable
Connectedness Property at any $\alpha,\beta \in \sper\ A$.
\end{conj}

Exactly the same proof which shows that the Connectedness Property
implies the Pierce-Birkhoff Conjecture applies verbatim to show that the
Definable Connectedness Property implies the Pierce-Birkhoff
conjecture for any ring $A$.
\medskip

One advantage of the Connectedness conjecture is that it is a
statement about $A$ (resp. about polynomials if $A=B$) which makes no
mention of piecewise polynomial functions.
\medskip

The Connectedness Conjecture is local in $\alpha$ and $\beta$. The purpose of
this paper is to associate to each pair $(\alpha,\beta)$ with
$ht(<\alpha,\beta>)=\dim A$ a natural number, called the complexity
of $(\alpha,\beta)$, and prove the Connectedness Conjecture in the simplest
case, according to this hierarchy, which is open : that of dimension 3 and
complexity 1.
\medskip

\begin{deft} Let $k$ be an ordered field. A $k$\textbf{-curvette} on
  $\sper\ A$ is a morphism of the form
$$
\alpha:A\to k\left[\left[t^\Gamma\right]\right],
$$
where $\Gamma$ is an ordered group. A $k$\textbf{-semi-curvette} is a
$k$-curvette $\alpha$ together with
  a choice of the sign data $\sgn\ x_1$,..., $\sgn\ x_r$, where
  $x_1,...,x_r$ are elements of $A$ whose $t$-adic values induce an
  $\F_2$-basis of $\Gamma/2\Gamma$.
\end{deft}

We explained in \cite{LMSS2} how to associate to a point $\alpha$ of $\sper\ A$
a $\bar k_\alpha$-semi-curvette. Conversely, given an
ordered field $k$, a $k$-semi-curvette $\alpha$ determines a prime
ideal $\mathfrak{p}_\alpha$ (the ideal of all the elements of $A$
which vanish identically on $\alpha$) and a total ordering on
$A/\mathfrak{p}_\alpha$ induced by the ordering of the ring
$k\left[\left[t^\Gamma\right]\right]$ of formal power series.\medskip

Below, we will often describe points in the real spectrum by
specifying the corresponding semi-curvettes.
\medskip

Let $(A,\mathfrak{m},R)$ be a regular local ring of dimension $n$ and
$\nu$ a valuation centered in $A$; let $\Phi=\nu(A
\setminus
\{0\})$; $\Phi$ is a well-ordered set. For an ordinal $\lambda$, let
$\gamma_\lambda$ be the element of $\Phi$ corresponding to $\lambda$.

\begin{deft}
A \textbf{system of approximate roots} of $\nu$ is a countable well-ordered set 
$\mathbf{Q}=\{Q_i\}_{i \in \Lambda}$, $Q_i \in A$,  minimal in the sense of
inclusion, satisfying
the following condition : \label{def-sar} for every $\nu$-ideal $I$
in $A$, we have
\begin{equation} \label{eq:valideal0S}
I=\left\lbrace\left. \prod_j Q_j^{\gamma_j}\ \right|\ \sum_j \gamma_j \nu(Q_j)
  \geq \nu(I) \right\rbrace A.
\end{equation}
 By definition, each $Q \in \mathbf{Q}$ comes equipped with additional data,
called the expression of $Q$ and denoted by $Ex(Q)$. The expression is a sum of
generalized monomials involving approximate roots
which  precede $Q$ in the given order. 
 
A system of approximate roots of $\nu$ up to $\gamma_\lambda$ is a
well-ordered set of elements of $A$ satisfying (\ref{eq:valideal0S})
only for $\nu$-ideals $I$ such that $\nu(I) < \gamma_\lambda$.

A finite product of the form $\mathbf{Q}^\eta=\prod\limits_jQ_j^{\eta_j}$ with
$\eta_j\in\mathbb N$ is called a \textbf{generalized monomial}. We order the
set of generalized monomials by the lexicographical order of the pairs
$(\nu(\mathbf{Q}^\eta),\eta)$ (cf. \cite{LMSS2}, below Definition 1.4).
\end{deft}

In  paragraph 1.2, Theorem 1.7 of \cite{LMSS2}, we constructed a system of
approximate roots up to some $\gamma$, $Q_i$, recursively in $i$. 
From now on, we fix this system of approximate roots once and for all.
\medskip

Let $u_1,\ldots,u_n$ be a regular system of parameters of $A$.

\begin{deft}
Let $i \in \N $ be a natural number, consider an approximate root
$(Q,Ex(Q))$. The notion of $Q$ being of complexity $i$ is defined as
follows. We say that $Q$ is an approximate root of complexity 0 if $Q
\in \{u_1,\ldots,u_n\}$. For $i>0$, we say that $Q$ is of complexity $i$ if
all the approximate roots appearing in $Ex(Q)$ are of complexity at most
$i-1$ and at least one approximate root appearing in $Ex(Q)$ is of
complexity precisely $i-1$.
\end{deft}

Fix $\alpha, \beta \in \sper\ A$ and consider the Connectedness conjecture for
this pair $(\alpha,\beta)$. Assume $\sqrt{<\alpha,\beta>}=\mathfrak{m}$. We
now define a natural number, called the complexity of $(\alpha,\beta)$.

\begin{deft}
 The complexity of  $(\alpha,\beta)$ is the smallest natural number $i$ such
that
every $\nu_\alpha$-ideal containing $<\alpha,\beta>$ is generated by
generalized monomials involving approximate roots of complexity at most $i$.
\end{deft}

In \cite{LMSS}, we proved the Connectedness conjecture for polynomial
rings of arbitrary dimension over a real closed field and pairs
$(\alpha,\beta)$ of complexity 0. Using Corollary \ref{prop-connexe} below,
based on \cite{And}, Chapter VII, 8.6, we can extend this result to
the case of excellent regular local rings $A$ of arbitrary dimension and pairs
$(\alpha,\beta)$ of complexity 0.

In this paper, we will assume that $R=\R$. In this case, $Ex(Q)$ is a binomial
in the approximate roots preceding $Q$ as we show below. The main result of
this paper is :
\begin{thm} \label{th:ppal}
 Let $(A,\mathfrak{m},\R)$ be an excellent 3-dimensional regular local ring
such that $\R \hookrightarrow A$. Let $\alpha, \beta \in \sper\ A$. Assume that
one of the following holds :

(1) $ht(<\alpha,\beta>) \leq 2$

(2) $ht(<\alpha,\beta>)=3$ and either $\{u_1,u_2,u_3\} \cap
<\alpha,\beta> \neq \emptyset$ or $(\alpha,\beta)$ is of complexity at most 1.

\noi Then the Connectedness Conjecture (and hence the Local Pierce-Birkhoff
Conjecture) holds for $(\alpha,\beta)$.
\end{thm}

Fix $\alpha,\beta \in \sper\ A$ and let
$\mathfrak{p}= \sqrt{<\alpha,\beta>}$. The case when $ht(\mathfrak{p})=1$ is
easy. The proof given in \cite{LMSS2} works verbatim in any dimension.

The present paper is organized as follows.

In \S\ref{strong2} we state a new conjecture, called the Strong Connectedness
Conjecture. We show that the Strong connectedness conjecture in dimension $n-1$
implies the Connectedness conjecture in dimension $n$ whenever
$ht(\mathfrak{p}) <\dim\ A$.

In \S\ref{ht2} we prove the Strong Connectedness Conjecture for arbitrary
regular local rings of dimension 2. We deduce the Connectedness Property
and the local Pierce-Birkhoff Conjecture for any ring $A$ (of any dimension)
and $\alpha,\beta \in \sper\ A$ such that $ht(<\alpha,\beta>)=2$ and
$A_{\sqrt{<\alpha,\beta>}}$ is regular.

\S\ref{sec:reel} is devoted to the study of graded algebras associated to
points of real spectra in the case when the residue field of our local ring is
$\R$. 

In \S\ref{monomial} we prove a comparison theorem between connected
components of a constructible subset $C \subset \sper\ A$ and those of the set
$\tilde{C} \subset \sper\ R[u_1,\ldots,u_n]_{(u_1,\ldots,u_n)}$
defined by the same formulae as $C$.

Finally, we describe some subsets of
$\sper\ A$, containing $\alpha$ and $\beta$, which will be later proved to be
connected, thus verifying the Connectedness Conjecture.

In \S\ref{ht3} we prove the Connectedness conjecture in the Case 2 of the
Theorem \ref{th:ppal}.

\section{The Strong Connectedness Conjecture}\label{strong2}

Let $\dim\ A=3$ and $ht\ \mathfrak{p}=2$. A natural idea would be
to apply the already known 2-dimensional connectedness conjecture
to the regular 2-dimensional local ring $A_\mathfrak{p}$. Then one would
construct a sequence of point blowings up $\tilde\pi:\tilde X_{\tilde
l}\rightarrow\sper\ A_\mathfrak{p}$ and a connected set in $\tilde X_{\tilde l}$
satisfying the conclusion of the conjecture. Finally, we would construct a
sequence $\pi:X_l\rightarrow\sper\ A$ of blowings up of points and smooth curves
whose restriction to the generic point of $V(\mathfrak p)$ is $\tilde\pi$.

The difficulty with this approach is that the 2-dimensional connectedness
conjecture cannot be applied directly. Indeed, let $g_1,\dots,g_s$ be as in the
connectedness conjecture and let $\Delta_\alpha\subset\Gamma_\alpha$ denote the greatest
isolated subgroup not containing $\nu_\alpha(\mathfrak p)$. \\ The hypothesis
$g_i\notin<\alpha,\beta>$ does not imply that
$g_i\notin<\alpha,\beta>A_{\mathfrak p}$: it may happen that
$\nu_\alpha(g_i)$ $< \nu_\alpha(<\alpha,\beta>)$,
$\nu_\alpha(g)-\nu_\alpha(\mathfrak p)\in\Delta_\alpha$ and so
$g_i\in<\alpha,\beta>A_{\mathfrak p}$, as we show by the example below.

\smallskip

\noindent\textbf{Example.} Let $\alpha,\beta$ be given by the curvettes
\begin{eqnarray}
x(t)&=&t^{(0,3)}\label{eq:xt3}\\
y(t)&=&t^{(0,4)}+bt^{(1,0)}\\
z(t)&=&t^{(0,5)}+ct^{(1,1)}\label{eq:zt5},
\end{eqnarray}
where $b\in\{b_\alpha,b_\beta\}\subset\R$ and
$c\in\{c_\alpha,c_\beta\}\subset\R$ and $t^{(0,1)}>0, t^{(1,0)}>0$. The
constants $b_\alpha,b_\beta,c_\alpha,c_\beta$ will be specified later. Let
$f_1=xz-y^2$, $f_2=x^3-yz$, $f_3=x^2y-z^2$; consider the ideal
$(f_1,f_2,f_3)$. The most general common specialization of $\alpha,\beta$ is given by the curvette
\begin{eqnarray}
x(t)&=&t^3\\
y(t)&=&t^4\\
z(t)&=&t^5,
\end{eqnarray}
$t>0$.
The corresponding point of $\sper\ A$ has support $(f_1,f_2,f_3)$, so $\mathfrak p=\sqrt{<\alpha,\beta>}=(f_1,f_2,f_3)$. Let $(x_\alpha(t),y_\alpha(t),z_\alpha(t))$ and $(x_\beta(t),y_\beta(t),z_\beta(t))$ be the curvettes defining $\alpha$ and $\beta$ as in (\ref{eq:xt3})--(\ref{eq:zt5}). Let us calculate $f_1(x_\alpha(t),y_\alpha(t),z_\alpha(t))$ and $f_1(x_\beta(t),y_\beta(t),z_\beta(t))$. In the notation of (\ref{eq:xt3})--(\ref{eq:zt5}) we have
\begin{eqnarray}
f_1(x(t),y(t),z(t))&=&(c-2b)t^{(1,4)}+\tilde f_1\\
f_2(x(t),y(t),z(t))&=&-(c+b)t^{(1,5)}+\tilde f_2\\
f_3(x(t),y(t),z(t))&=&(b-2c)t^{(1,6)}+\tilde f_3,
\end{eqnarray}
where $\tilde f_i$ stands for higher order terms with respect to the
$t$-adic valuation. Choose $b_\alpha, b_\beta, c_\alpha,$ $c_\beta$ so that
none of $f_1,f_2,f_3$ change sign between $\alpha$ and $\beta$. The
smallest $\nu_\alpha$ value of an element which changes sign between
$\alpha$ and $\beta$ is $(1,4)+(0,4)=(1,5)+(0,3)=(1,8)$, so
$\nu_\alpha(<\alpha,\beta>)=(1,8)$. Thus we have
$f_i\notin<\alpha,\beta>$, but $f_i\in<\alpha,\beta>A_{\mathfrak p}$,
as desired.
\medskip

Thus we are naturally led to formulate a stronger version of the Connectedness
Conjecture, one which has exactly the same conclusion but with
somewhat weakened hypotheses. This phenomenon occurs in all dimensions, as we
now explain.
\medskip

\begin{deft} \textbf{Strong Connectedness Property}
Let $\Sigma$ be a ring, $\alpha,\beta \in \sper\ \Sigma$, having a common
specialization $\xi$.
We say that $\Sigma$ has the Strong Connectedness
Property at $\alpha, \beta$ if given any $g_1,\ldots,g_s \in \Sigma \setminus
(\mathfrak{p}_\alpha \cup \mathfrak{p}_\beta)$ such
  that for all $j \in \{1,\ldots,s\}$,
\begin{equation} \label{eqn:sdcc}
 \nu_\alpha(g_i) \leq
  \nu_\alpha(<\alpha,\beta>),\ \ \nu_\beta(g_i) \leq
  \nu_\beta(<\alpha,\beta>)
\end{equation}
 and such that no $g_i$ changes sign between $\alpha$ and $\beta$, the points
$\alpha$ and $\beta$ belong to the same connected component of $\sper\
\Sigma \setminus \{g_1 \cdots g_s =0\}$.
\end{deft}

\begin{conj} \textbf{Strong Connectedness Conjecture} \label{strong}
Let $\Sigma$ be a regular ring. Then $\Sigma$ has the Strong Connectedness
Property at any pair of points $\alpha,\beta \in \sper\ \Sigma$
having a common specialization.
\end{conj}

 Let $A$ be a ring and $\alpha, \beta \in \sper\ A$. Let $\mathfrak{p} =
\sqrt{<\alpha,\beta>}$, let $\alpha_0$ be the pre-image of $\alpha$ under the
natural inclusion $\sigma : \sper\ A_{\mathfrak p} \hookrightarrow \sper\
A$ and similarly for $\beta_0$.

\begin{thm} \label{thm:localsdcc}
If $\sper\ A_{\mathfrak{p}}$ has the Strong Connectedness property at
$\alpha_0,\beta_0$, then $A$ satisfies the Connectedness Conjecture at
$\alpha, \beta$.
\end{thm}

\noi Proof : Let $g_1,\dots,g_s\in A$ be the elements appearing in the statement
of the Connectedness Conjecture. Renumbering the $g_i$, if necessary,
we may assume that $g_1,\dots,g_l\notin<\alpha_0,\beta_0>$ and
$g_{l+1},\dots,g_s\in<\alpha_0,\beta_0>$. The condition
$g_{l+1},\dots,g_s\in<\alpha_0,\beta_0>$ implies that, for
$i \in \{l+1,\ldots,s\}$,
$\nu_{\alpha_0}(g_i)=\nu_{\alpha_0}(<\alpha_0,\beta_0>)$.

By hypothesis, there exists a connected set $C_0 \subset
\sper\ A_{\mathfrak{p}}$, $\alpha_0,\beta_0 \in C_0$
such that $C_0 \subset \{g_1 \cdots g_s \neq 0\}$. Then $\sigma(C_0)$
satisfies the conclusion of the Connectedness Conjecture for
$A,\alpha,\beta,g_1,\ldots,g_s$.
$\Box$ \medskip

In the next section we will use Zariski's
theory of complete ideals to prove the Strong Connectedness Conjecture in
dimension 2, and hence also the Connectedness conjecture in dimension 3, when
$ht(\mathfrak{p})=2$.

\section{The case when the height of $\mathfrak{p}$ is 2}\label{ht2}

\begin{thm} \label{thm:sdcc}
Conjecture \ref{strong} is true when $\Sigma$ is of dimension 2.
\end{thm}

\noi Proof : If one of $\alpha$, $\beta$ is a specialization of the other, the
result is trivially true, because the connected component of $\sper\
\Sigma \setminus \{g_1 \cdots g_s  =0\}$ containing the more general point among
$\alpha$ and $\beta$ satisfies the conclusion of the conjecture. From now on we
shall assume that none of $\alpha$ and $\beta$ is a specialization of the other.
\smallskip

Let $z$ be a new variable. We will say that a point $\eta \in \sper\ k[z]$ is
closed if $\{\eta\}=\sur{\{\eta\}}$.

Consider a point $\alpha \in \sper\ \Sigma$, $\dim\ \Sigma =2$. Let $\xi$ be
the most special specialization of $\alpha$. Assume that ht$(\mathfrak{p}_\xi) =
2$ and $\alpha \neq \xi$. Let $(x,y)$ be a regular
system of parameters of $\Sigma_{\mathfrak{p}_\xi}$ and let $k$ be the residue
field $k = \frac{\Sigma}{\mathfrak{p}_\xi}$.
Let $\rho: X \to \sper\ \Sigma$ be the blowing up of $\sper\ \Sigma$ along
$(x,y)$. Let $\alpha'$ be the strict transform of $\alpha$ in $X$ (see
\cite{LMSS2}, Definitions 3.19 and 3.20). If $\nu_\alpha(y) \geq
\nu_\alpha(x)$ then $\alpha' \in \sper\ \Sigma[\frac{y}{x}]$. Consider the
homomorphism $\Sigma[\frac{y}{x}] \to k[z]$ which maps $\frac{y}{x}$ to $z$ and
elements of $\Sigma$ to their image in $k$. In this way, we identify $\sper\
k[z]$ with $\sper\ \Sigma[\frac{y}{x}]\cap \rho^{-1}(\xi)$.

\begin{deft} \label{slope}
The slope of $\alpha$, denoted by $sl(\alpha)$, is the following element of
$\sper\ k[z] \cup \{\infty\}$

- if $\nu_\alpha(x) > \nu_\alpha(y)$, $sl(\alpha) := \infty$;

- if  $\nu_\alpha(x) \leq \nu_\alpha(y)$, $sl(\alpha)$ is the most special
specialization of $\alpha'$ in $\sper\ \Sigma[\frac{y}{x}]$.

Let $\alpha, \beta \in \sper\ \Sigma$ be the two points centered at $\xi$ and
having the same slope. We say that $\alpha$ and $\beta$ point in the same
direction if sgn$(x(\alpha)) = $ sgn$(x(\beta))$ when $sl(\alpha) \neq
\infty$ (resp.  sgn$(y(\alpha)) = $ sgn$(y(\beta))$ when $sl(\alpha) =
\infty$). Otherwise we say that $\alpha$ and $\beta$ point in different
direction.
\end{deft}

\noi Examples : Let $\Sigma= \Q[x,y]$.

1. Let $\alpha$ be the point of $\sper\ \Sigma$ given by the following
semi-curvette $\Q[x,y] \hookrightarrow \Q(\pi)[[t]]$ such that $x \mapsto t,\ y
\mapsto \pi t$. Then $\xi$ is the closed point with support $(x,y)$ and the
slope of $\alpha$ is the point of $\sper\ \Q[z]$ such that for any rational
number $p/q$ we have $z > p/q \iff \pi > p/q$.

2. Let $\alpha$ be a point of $\sper\ \Sigma$ such that
$\nu_\alpha(x)=\nu_\alpha(y) >0$, $\nu_\alpha(y^2-2 x^2) > 2 \nu_\alpha(x)$.
Then $\xi$ is the closed point with support $(x,y)$ and the
slope of $\alpha$ is the point of $\sper\ \Q[z]$ with support $(z^2-2)$.
\medskip

\begin{rek}
 In the situation of Definition \ref{slope}, assume that $sl(\alpha) \neq
\infty$. Let $k[z](sl(\alpha))$ be the field of fractions of
$\frac{k[z]}{\mathfrak{p}_{sl(\alpha)}}$. We can naturally identify
$k[z](sl(\alpha))$ with
the ordered sub-field of $k_\alpha$ generated over $k$ by the image of
$\frac{y}{x}$. The field $k[z](sl(\alpha))$ is a simple extension of $k$ which
can be algebraic as in the Example 2 above, or transcendental as in the Example
1.
\end{rek}

\begin{deft}
 Let $f \in \sper\ \Sigma$. We say that $f=0$ is tangent to $\alpha$ if $
\nu_\alpha(f) > \nu_\alpha(\mathfrak{p}_\xi)$.
\end{deft}

First assume that $\alpha$ and $\beta$ have the same tangent, and that they
are facing in different directions along that tangent. Then $<\alpha,\beta>
=\mathfrak{p}_\xi$. We want to show that, for all $i$, $g_i \notin
\mathfrak{p}_\xi$. Assume that $g_i \in \mathfrak{p}_\xi$. Write $g_i=ax+by +
\tilde{g_i}$ where $a,b \in \Sigma$ and $\tilde{g_i} \in (x,y)^2$. We may
assume that the common slope to $\alpha$ and $\beta$ is not $\infty$. Then
\begin{eqnarray} \label{eqn:slope}
\nu_\alpha(g_i) &= &\nu_\alpha(\mathfrak{p}_\xi)=\nu_\alpha(x) \leq
\nu_\alpha(y) \\ \label{eqn:slope2}
\nu_\beta(g_i) &= &\nu_\beta(\mathfrak{p}_\xi)=\nu_\beta(x)
\leq \nu_\beta(y).
\end{eqnarray}
Hence either $a \notin \mathfrak{p}_\xi$ or $(\nu_\alpha(x)
= \nu_\alpha(y))$ and $b \notin \mathfrak{p}_\xi$. In particular,
$sgn_\alpha(g_i) = sgn _\alpha(ax+by)$ and similarly for $sgn_\beta$.

Let $k[z](sl(\alpha))$ be as in the previous remark.
By (\ref{eqn:slope}) and (\ref{eqn:slope2}), the natural image of
$a+b\frac{y}{x}$ in $k[z](sl(\alpha))$is non zero. Since $\alpha$ and
$\beta$ have the same
slope, they induce the same order on $k[z](sl(\alpha))$. Hence
$a+b\frac{y}{x}$ does not
change sign between $\alpha$ and $\beta$, so $x(a+b\frac{y}{x})$ changes sign
between $\alpha$ and $\beta$, which is a contradiction. Hence $g_i \notin
\mathfrak{p}_\xi$. Then a small connected neighbourhood $U$ (small enough so
that
$\{g_1\cdots g_s=0\}\cap U = \emptyset$) of $\xi$ satisfies the conclusion of
Conjecture \ref{strong}. This proves the Theorem in the special case when
$\alpha$ and $\beta$ have the same slope but point in different directions.
\smallskip

\noi From now on assume that if $\alpha$ and $\beta$ have the same slope,
they point in the same direction.
\smallskip

Let $\pi: X' \to X = \sper\ A$ the shortest sequence of blowings up such that
the strict transforms $\alpha'$ and $\beta'$ of
$\alpha$ and $\beta$ have the same specialization $\xi'$ with
$ht(\mathfrak{p}_{\xi'}) = 2$ and distinct slopes
(see \cite{LMSS2}, by iterating Proposition 3.31). Note
that, if $g'_i$ denotes the strict transform of $g_i$, then the $g'_i$ such that
$g_i'(\xi) \neq 0$ play no role and if $g'_i(\xi)=0$, by (\ref{eqn:sdcc}),
$\{g'_i=0\}$ cannot be tangent to $\alpha'$ or $\beta'$ or to the last
exceptional divisor if it exists. Let
$\mathcal{O}_{X',\xi'}$ be the local ring of $X'$ at $\xi'$ and let
$x',y'$ be a regular system of parameters such that $\{x'=0\}$ is the
last exceptional divisor if it exists and $\{y'=0\}$ the second one if
any. In the case we had not to blow up, we take an $x'$ such that
$\{x'=0\}$ is not tangent to $\alpha',\beta'$ or any of $\{g'_i=0\}$ and such
that $x'(\alpha')>0$ and $x'(\beta') >0$.
Note that $x'$ does not change sign between $\alpha$ and $\beta$ (otherwise the
blowing up sequence $\pi$ would have stopped at an earlier stage). Replacing
$x'$ by $-x'$ if necessary, we may assume that $x'(\alpha') > 0$, $x'(\beta')
>0$.
\smallskip

Let us introduce the following total ordering on the set
$\{g'_1,\ldots,g'_s\}$. Write each $g'_j$ as a formal power
series in the formal completion $\mathcal{O}_{X',\xi'} \to k'[[x',y']]$
as $$g'_j=y'+\sum_{i=1}^\infty c_{ij} {x'}^i \mbox{ with } c_{ij} \in k'.$$
This is possible because of the choice of $x',y'$, the non tangency of
$\{g'_j=0\}$ with $\alpha', \beta'$ and the last exceptional divisor. We compare
$g'_j$ and $g'_\ell$ by comparing the monomials in
lexicographic ordering. Namely, we take the smallest $i$ such that
$c_{ij} \neq c_{i\ell}$ and we say that $j \prec \ell$ if $c_{ij} <
c_{i\ell}$. Without loss of generality, we may assume that
$g'_1(\alpha') >0,\ldots,g'_\ell(\alpha')>0$,
$g'_{\ell+1}(\alpha')<0,\ldots,g'_s(\alpha')<0$ and also that
$1 \prec \ldots \prec \ell$, $\ell+1 \prec \ldots \prec s$.

\begin{lem} \label{lem-inclus}
Let $j,q \in \{1,\ldots,\ell\}$, $j \prec q$. Then $\{g'_j>0,x'>0\}
\subset \{g'_q>0,x'>0\}$.

Let $j,q \in \{\ell+1,\ldots,s\}$, $j \prec q$. Then $\{g'_q>0,x'>0\}
\subset \{g'_j>0,x'>0\}$.
\end{lem}

\noi Proof : In the first case, we have to prove that $g'_q(\delta) > 0
\Ra g'_j(\delta)>0$. Write $g'_j = y'+ c_{1j}x'+ \cdots+ c_{ij}{x'}^i +
{x'}^{i+1}( \cdots)$ and $g'_q=  y'+ c_{1q}x'+ \cdots+ c_{iq}{x'}^i +
{x'}^{i+1}( \cdots)$ with $c_{kj}=c_{kq}$ for $k=1,\ldots,i-1$ and
$c_{ij} < c_{iq}$. We have $g'_q-g'_j= (c_{iq}-c_{ij}){x'}^iu$ where $u$
is a positive unit of $k'[[x',y']]$. So $g'_q-g'_j=d{x'}^i$ in
$\mathcal{O}_{X',\xi'}$ with $d \in \mathcal{O}_{X',\xi'} \setminus
\mathfrak{m}_{X',\xi'}$ such that $d = c_{iq}-c_{ij} \mod
\mathfrak{m}_{X',\xi'}$,
in particular $d(\delta)>0$.

And the same with the second inclusion. $\Box$
\medskip

\begin{lem}
We have $c_{11} > c_{1s}$.
\end{lem}

\noi Proof : Note that we have
\begin{eqnarray}
g'_1(\alpha') = y'(\alpha') + c_{11}x'(\alpha') + x'(\alpha')^2 h_1 > 0
\label{ineq1} \\
g'_1(\beta') = y'(\beta') + c_{11}x'(\beta') + x'(\beta')^2 h_1 > 0
\label{ineq2} \\
g'_s(\alpha') = y'(\alpha') + c_{1s}x'(\alpha') + x'(\alpha')^2 h_s < 0
\label{ineq3}\\
g'_s(\beta') = y'(\beta') + c_{1s}x'(\beta') + x'(\beta')^2 h_s < 0
 \end{eqnarray}
where $h_1, h_s \in k'[[x',y']]$.


Write $\alpha'$ as curvette :
\begin {eqnarray*}
x'(t) &=& t^{\nu_\alpha(x')} + \cdots \\
y'(t) &=& b_\alpha t^{\nu_\alpha(x')} + \cdots
\end{eqnarray*}
where $b_\alpha$ is the natural image of $\frac{y'}{x'}$ in $k_\alpha$.

Then $g'_1(\alpha')>0 \Lra y'+c_{11}x'+{x'}^2h_1 > 0 \Lra y'(t) + c_{11} x'(t)
+ x'(t)^2h_1(t) = (b_\alpha+c_{11})t^{\nu_\alpha(x')} + \cdots > 0$ in
$\sur{k_\alpha}[[t^{\Gamma_\alpha}]]$, so $b_\alpha +c_{11} \geq 0$ in
$k_\alpha$. By the same
arguments, we have : $b_\alpha+c_{1s} \leq 0$ in $k_\alpha$, $b_\beta + c_{11}
\geq 0$ and $b_\beta+c_{1s} \leq 0$ in $k_\beta$. If we had $b_\alpha + c_{11}
=0$ and $b_\beta+c_{11}=0$, then $b_\alpha=b_\beta = -c_{11} \in k'$. Hence
$k[z](sl(\alpha')) = k[z](sl(\beta'))$, which contradicts the fact that
$\alpha'$ and $\beta'$ have different slopes. Thus,
at least one of the inequalities (\ref{ineq1}) and (\ref{ineq2}) is
strict, say $b_\alpha + c_{11} > 0$ for instance. Together with the
inequality (\ref{ineq3}), this implies that $c_{11} > c_{1s}$. $\Box$
\medskip

Let $$C'= \{g_1'>0,g_s'<0\}.$$ By
definition $C'$ contains $\alpha'$ and $\beta'$, so is non empty.

\begin{cor}
For any $\delta \in C'$, $x'(\delta) >0$.
\end{cor}

\noi Proof : It's a straightforward consequence of the preceding Lemma and
proof. Of course, $g_1'(\delta) > 0 \Ra x'(\delta)(b_\delta+c_{11}) \geq 0$
and $g'_s(\delta) <0 \Ra x'(\delta)(b_\delta +c_{1s}) \leq 0$ where $b_\delta$
is defined in a similar way as $b_\alpha$ or $b_\beta$. So
$x'(\delta)(c_{11}-c_{1s}) \geq 0$, which proves the result.
$\Box$
\medskip

By Lemma \ref{lem-inclus}, for any $\delta \in C$, we have
$g'_1(\delta) >0,\ldots,g'_\ell(\delta)>0$,
$g'_{\ell+1}(\delta)<0,\ldots,g'_s(\delta)<0$. So, finally, being a
quadrant in $\sper\ \mathcal{O}_{X',\xi'}$, if $\mathcal{O}_{X',\xi'}$ is an
excellent ring, $C'$ is connected by
(\cite{LMSS2}, Theorem 3.35). So the image $C$ of $C'$ in $X$ satisfies the
conclusion of the Strong Connectedness Conjecture. $\Box$

\begin{cor}
Let $A$ be a ring, $\alpha,\beta \in \sper\ A$, $\mathfrak{p}=
\sqrt{<\alpha,\beta}$. Assume that the local ring $A_\mathfrak{p}$ is
excellent regular
of dimension at most 2. Then $A$ has the Connectedness Property at
$\alpha,\beta$ and hence satisfies the local Pierce-Birkhoff conjecture at
$\alpha,\beta$.
\end{cor}

This follows immediately from Theorem \ref{thm:localsdcc} and Theorem
\ref{thm:sdcc}.

\begin{rek}
 All the results of this section remain true, with the same proofs if we drop
the excellence hypothesis on $A$, but replace the Strong Definable
Connectedness Conjecture by Definable Strong Connectedness Conjecture.
\end{rek}

\section{Graded algebras in the case of residue field $\R$} \label{sec:reel}
\begin{thm}\label{th:reel}
  Let $(\Sigma,\mathfrak{m})$ be a local ring with residue field $\R$. Let
$\alpha \in \sper\ \Sigma$ such that $\nu_\alpha$ is centered at
$\mathfrak{m}$. For every $\gamma \in \Gamma_\alpha$, we have
$\frac{P_\gamma}{P_{\gamma_+}} \cong \R$. 
\end{thm}

\noi Proof : We have non-canonical inclusions $$\R \subset
\frac{P_\gamma}{P_{\gamma_+}} \subset k_\alpha.$$ Thus it is sufficient to prove
that $\R \equiv k_\alpha$. Take an element $\sur{b} \in k_\alpha$ and let $b$
be a representative of $\sur {b}$ in $R_\alpha$. By definitions, there exists
$a \in \Sigma[\alpha]$ such that $|b| \leq a$. Since $\nu_\alpha$ is centered at
$\mathfrak{m}$, we may take $a$ to be a unit of $\Sigma[\alpha]$. Let $\sur{a}$
be the image of $a$ in $\frac{\Sigma[\alpha]}{\mathfrak{m}\Sigma[\alpha]} \cong
\R$. Then $|\sur{b}| \leq \sur{a}$. Hence $\sur{b} \in \R$ as desired. $\Box$ 

\begin{cor}\label{unidiminsionnel}
Assume in addition that $\Sigma$ is regular and $\R \subset \Sigma$. Let $Q$ be
an
approximate root for $\alpha$, then $Ex(Q)$ is a difference of two generalized
monomials in the approximate roots preceding $Q$. 
\end{cor}

\noi Proof : By the construction of approximate roots (\cite{LMSS2}, section
1.2), $Ex(Q)$ comes from a certain $\R$-linear dependence relation among
generalized monomials in approximate roots, preceding $Q$, of the same value.
According to Theorem \ref{th:reel}, any two such monomials are $\R$-multiples of
each other. Now the result follows from the construction of approximate roots.
$\Box$

\section{Connectedness properties}\label{monomial}
\begin{thm} \label{th:connex}
Let $(A,\mathfrak{m},R)$ be an excellent regular local ring such that $R
\subset A$. Let $(u_1,\ldots$ $,u_n)$ be a regular system of parameters
 of $A$. Let $C \subset \sper\ A$ be a constructible set such that all the
elements of
$A$ appearing in formulae defining $C$ belong to
$R[u_1,\ldots,u_n]_{(u_1,\ldots,u_n)}$. Let $\tilde{C} \subset \sper\
R[u_1,\ldots,u_n]_{(u_1,\ldots,u_n)}$ denote the constructible set defined by
the same formulae as
$C$. Let $$U = \{ \delta \in \sper\ A\ |\ \delta \text{ is centered at }
\mathfrak{m}\}$$ $$\tilde{U} = \{ \delta \in \sper\
R[u_1,\ldots,u_n]_{(u_1,\ldots,u_n)} \ |\ \delta \text{ is centered at }
(u_1,\ldots,u_n) \}.$$ Then the natural map $\sper\ A \to \sper\
R[u_1,\ldots,u_n]_{(u_1,\ldots,u_n)}$ induces a bijection between the set of
connected components of $C\cap U$ and the set of connected components of
$\tilde{C} \cap \tilde{U}$.
\end{thm}

\noi Proof : Consider the
following natural ring
homomorphisms $$\xymatrix{R[u_1,\ldots,u_n]_{(u_1,\ldots,u_n)}
\ar[r]^-{\sigma_0} & A \ar[r]^-\sigma & R[[u_1,\ldots,u_n]]}.$$

The theorem follows from (\cite{And}, chap. VII, Proposition 8.6) applied to
the rings $A$ and $R[u_1,\ldots,u_n]_{(u_1,\ldots,u_n)}$.
$\Box$

%
%
%
%
%
%
%
%


\begin{cor} \label{prop-connexe}
Let $(A,\mathfrak{m},R)$ be an excellent regular local ring such that $R
\subset A$. Let $(u_1,\ldots,u_n)$ be a regular system of parameters
 of $A$. Fix a subset $J \subset \{1,\dots,n\}$ and the point $\xi \in \sper\ A$
such that $\mathfrak{p}_\xi=\mathfrak{m}$.
Let $U$ denote the subset of $\sper\ A$ consisting of generizations of $\xi$.
Let $C$ denote the subset of $U$ defined by specifying $\sgn\ u_q$ (which can be
either strictly positive on all of $C$ or strictly negative on all of $C$) for
$q \in J$ and by imposing, in addition, finitely many monomial inequalities of
the form
\begin{equation}
\left| d_i \sous{u}^{\lambda_i} \right| \ge \left| \sous{u}^{\theta_i}
\right|,\quad 1 \le
i \le M  \label{eq:monineqC1}
\end{equation}
where $d_i \in R\setminus\{0\}$, $\lambda_i,\theta_i\in\N^n$ and $u_q$ may
appear only on the right hand side of the inequalities (\ref{eq:monineqC1}) for
$q \notin J$. Then $C$ is connected.
\end{cor}

\noi Proof : Write $\lambda_i=(\lambda_{1i},\ldots, \lambda_{ni})$ 
and similarly for $\theta_i$.
It is sufficient to prove that any two elements of
$C$ belong to the same connected component of
$C$.

Consider the natural homomorphism
\begin{equation}
A \to \hat{A}= R[[u_1,\ldots,u_n]].
\end{equation}

Let $\hat{\xi}$ denote the point of $\sper\ \hat{A}$ with support
$\mathfrak{m}\hat{A}$.

Following (\cite{And}, chap. VII, proposition 8.6), $C$ is
connected if
and only if $$\hat{C} = \{\delta \in \sper\
R[[u_1,\dots,u_n]] \ | \ u_j(\delta)>0, \ j \in J, \left| d_i
\sous{u}^{\lambda_i} \right| \ge \left| \sous{u}^{\theta_i}
\right|,\ 1\le i \le M, \hat{\xi} \in \sur{\{\delta \}} \}$$
is connected (this is where we are using the fact that $A$ is
excellent). So it suffices to prove that $\hat{C}$ is
connected.

By the preceding Theorem, $\hat{C}$ is connected if and only if the set
\begin{eqnarray*}
C_\dag &=& \{\delta \in \sper\
R[u_1,\dots,u_n]_{(u_1,\dots,u_n)} \ | \
u_j(\delta)> 0, \ j \in J,\\ & &\left| d_i \sous{u}^{\lambda_i} \right| \ge
\left| \sous{u}^{\theta_i} \right|,\ 1\le i \le M,\delta \mbox{ is centered at
}(u_1,\dots,u_n) \}
\end{eqnarray*}
is connected.

Define
\begin{eqnarray*}
C_0 &=& \{\delta \in \sper\ R[u_1,\dots,u_n] \ | \
u_j(\delta)> 0, \ j \in J,\\ & &\left| d_i \sous{u}^{\lambda_i} \right| \ge
\left| \sous{u}^{\theta_i} \right|,\ 1\le i \le M,\delta \mbox{ is centered at
}(u_1,\dots,u_n) \}
\end{eqnarray*}
and
\begin{eqnarray*}
C_{loc} &=& \{\delta \in \sper\ R[u_1,\dots,u_n]_{\prod_{j \in J} u_j} \ | \
u_j(\delta)> 0, \ j \in J,\\ & &\left| d_i \sous{u}^{\lambda_i} \right| \ge
\left| \sous{u}^{\theta_i} \right|,\ 1\le i \le M,\delta \mbox{ is centered at
}(u_1,\dots,u_n) \}
\end{eqnarray*}

\noi The natural maps $\phi:R[u_1,\ldots,u_n] \to
R[u_1,\ldots,u_n]_{(u_1,\ldots,u_n)}$ and $\psi : R[u_1,\ldots,u_n] \to
R[u_1,\ldots,u_n]_{\prod_{j\in J} u_j}$ induce homeomorphisms
$\phi|_{C_0}:C_0 \cong C_{loc}$ and $\psi|_{C_0}:C_0 \cong C_\dag$.

So it suffices to prove that $C_{loc}$ is connected. But
$$
C_{loc} =
\bigcap_{N \in \N} C_N
$$ where

\begin{eqnarray*}
C_N &=& \{\delta \in \sper\ R[u_1,\dots,u_n]_{\prod_{j \in J} u_j} \
| \ \frac{1}{N} \geq u_j(\delta) \geq  0, \ j \in J,\\ & & \left| d_i
\sous{u}^{\lambda_i} \right| \ge \left| \sous{u}^{\theta_i}
\right|,\ 1\le i \le M \}.
\end{eqnarray*}

By Lemma 4.1 of \cite{LMSS}, each $C_N$ is a non empty closed
connected subset of \\ $\sper\
R[u_1,\dots,u_n]_{\prod_{j\in J} u_j}$, hence $C_{loc}$ is connected
by (\cite{LMSS}, lemma 7.1). $\Box$
\medskip

\begin{rek} \label{rk:unit}
 Keep the hypothesis of Corollary \ref{prop-connexe}. Consider a set
$\tilde{C}$ defined by inequalities 
\begin{equation}
\left|\tilde{d}_i\sous{u}^{\lambda_i} \right| \geq \left| \sous{u}^{\theta_i}
\right|, \ 1\leq i\leq M, \ \tilde{d}_i \in A \setminus \mathfrak{m} 
\end{equation}
and the same sign conditions as $C$. For each $i$, $1\leq i\leq M$, take $d'_i
\in R$ such that $\left| \tilde{d}_i(\xi)\right| > \left|d'_i \right|$. Let
$\tilde{\tilde{C}}  \subset U$ be defined by
$\left|d'_i\sous{u}^{\lambda_i} \right| \geq \left| \sous{u}^{\theta_i}
\right|, \ 1\leq i\leq M$ and the same sign conditions as before. Then
$\tilde{\tilde{C}}$ is connected and $\tilde{\tilde{C}} \subset C$. 
\end{rek}

Assume that $A$ is of dimension 3 and has residue field $\R$. Let $\alpha, \beta
\in \sper\ A$ and suppose $ht(<\alpha,\beta>)$ $=3$. 

Let $\nu_{\alpha 0}$ be the monomial valuation defined by 
\begin{eqnarray}
\nu_{\alpha 0}(u_1) = \nu_\alpha(u_1) \\ 
\nu_{\alpha 0}(u_2) = \nu_\alpha(u_2) \\ 
\nu_{\alpha 0}(u_3) = \nu_\alpha(u_3).                                  
\end{eqnarray}
In other words, for a polynomial $f= \sum_{\gamma \in \N^3} c_\gamma
\sous{u}^\gamma$, we have $\nu_{\alpha 0}(f) =\min_\gamma
\{\nu_\alpha(\sous{u}^\gamma) \ |\ c_\gamma \neq 0\}$. 

\begin{deft}
 An approximate root $Q$ for $\alpha$, of complexity at most one, is said to be
relevant for
$(\alpha,\beta)$ if either $Q \in \{u_1,u_2,u_3\}$ or $\nu_{\alpha 0}(Q) <
\nu_\alpha(<\alpha,\beta>)$.
\end{deft}

Note that if $Q$ is relevant for $(\alpha,\beta)$, then
$Q$ is an approximate root for $\beta$. If, in addition $Q \notin
\{u_1,u_2,u_3\}$, then $\nu_{\beta 0}(Q) < \nu_\beta(<\alpha,\beta>)$.
\medskip

Let $\{Q_i\}_{4\leq i \leq \ell}$ where $\ell \in \{3,4,\ldots,r\}$ denote the
set of relevant approximate roots of complexity 1 (the case $\ell=3$ means that
the $(\alpha,\beta)$ has complexity 0). 
\medskip

Let
$g_1,\ldots,g_s \in A$ be as in the statement of the Connectedness Conjecture.
Let
\begin{equation} \label{eqn-gi}
g_i = \mathbf{Q}^{\delta_i} +
\sum_{j=1}^{N_i}c_{ji}\mathbf{Q}^{\delta_{ji}},\ i\in \{1,\ldots,s\}
\end{equation}
be the standard form of $g_i$ common to $\alpha$ and $\beta$ of level
$\nu_\alpha(<\alpha,\beta>)$ (see \cite{LMSS2},
\S 1.3); by definition then  $\mathbf{Q}^{\delta_i}$,
$\mathbf{Q}^{\delta_{ji}}$ are
generalized monomials in the relevant approximate roots and  
$\nu_\alpha(\mathbf{Q}^{\delta_i})<\nu_\alpha(\mathbf{Q}^{\delta_{ji} } )$ ,
$j\in\{1,\dots,N_i\}$. The fact that there is only one dominant monomial
$\mathbf{Q}^{\delta_i}$ is due to Theorem \ref{th:reel}.
\medskip

\noi 1. Let \begin{equation} \label{ensCC}
C = \left\lbrace \delta \in \sper A
    \ \left|\begin{array}{ll}
 \delta \text{ is centered at } (x,y,z) \\
\nu_\delta( \mathbf{Q}^{\delta_i}) <
        \nu_\delta(\mathbf{Q}^{\delta_{ji}})\ \forall i \in \{1,\ldots,s\}, \
\forall
 j \in \{1,\ldots,N_i\}
    \\
     sgn_\delta(Q_q) = sgn_\alpha(Q_q)  \mbox{ for all } Q_q
     \mbox{ appearing in } \mathbf{Q}^{\delta_i} 
 \end{array} \right. \right\rbrace.
\end{equation}

\noi 2. Let $C'$ defined by the set of all $\delta$, centered at
$(x,y,z)$, satisfying the inequalities
\begin{equation} \label{eq:dominant}
\left|\mathbf{Q}^{\delta_i}(\delta) \right| >
N_i|\mathbf{Q}^{\delta_{ji}}(\delta)| \ \forall i \in \{1,\ldots,s\}, \ \forall
j \in \{1,\ldots,N_i\}
\end{equation}
and the sign conditions appearing in (\ref{ensCC}).

\begin{rek}\label{Cworks} 1. We have $\alpha, \beta \in C$.

2. $C \cap \{g_1\cdots g_s = 0\} = \emptyset$. Indeed,
inequalities (\ref{ensCC}) imply that, for every $\delta \in
C$, $g_i(\delta)$ has the same sign as $\mathbf{Q}^{\delta_i}(\delta)$; in
particular, none of the $g_i$ vanish on $C$.

3. To prove the Connectedness Conjecture it is sufficient to prove that
$\alpha$ and $\beta$ lie in the same connected component of $C$.
\end{rek}
All the preceding remarks apply to $C'$.

\section{The case when $ht(\mathfrak{p})=3$}\label{ht3}

Let $(A,\mathfrak{m})$ be an regular local ring of dimension 3 with
residue field $\R$ contained in $A$. Let $u_1,u_2,u_3$ be a regular system of
parameters of $A$ such that $\nu_\alpha(u_1) \leq \nu_\alpha(u_2) \leq
\nu_\alpha(u_3)$. Let $\alpha,\beta \in \sper\ A$ centered at $\mathfrak{m}$. 
In this case, the approximate roots of complexity 1 are binomials in
$u_1,u_2,u_3$ (Corollary \ref{unidiminsionnel}).

\begin{lem}
Every valuation $\nu$ admits at most three approximate roots of complexity one.
\end{lem}

\noi Proof : An approximate root of complexity 1 is an irreducible binomial
$\omega_1-\omega_2$ having the property that $\nu(\omega_1-\omega_2) >
\nu(\omega_1)=\nu(\omega_2)$. We now prove that there are only three
possible types of approximate roots of complexity 1 (up to exchanging the 2
monomials in order to respect the monomial ordering defined in \cite{LMSS2},
\S 1.2, after Definition 1.4), that means :
\begin{eqnarray}
u_2^{\beta_1} u_3^{\gamma_1}-\lambda_1 u_1^{\alpha_1} \label{eqn-approot1} \\
u_2^{\beta_2} -\lambda_2u_1^{\alpha_2} u_3^{\gamma_2} \label{eqn-approot2}\\
u_3^{\gamma_3} -\lambda_3u_1^{\alpha_3} u_2^{\beta_3} \label{eqn-approot3}
\end{eqnarray}
$\lambda_1, \lambda_2,\lambda_3 \in \R$, with  $\alpha_1, \beta_2, \gamma_3$ the
smallest exponents possible.

By definition of approximate roots, the initial monomial of one
approximate root cannot be divisible by the initial monomial of
another. Therefore, there is at most one approximate root of each of
the forms  $u_2^{\beta_2} -\lambda_2 u_1^{\alpha_2} u_3^{\gamma_2}$ and
$u_3^{\gamma_3} -\lambda_3 u_1^{\alpha_3} u_2^{\beta_3}$.

We claim that there is also at most one
of the form $u_2^{\beta_1} u_3^{\gamma_1}-\lambda_1 u_1^{\alpha_1}$. Indeed
suppose there were another one
$u_2^{\beta'_1}u_3^{\gamma'_1}-\lambda_1' u_1^{\alpha'_1}$; then necessarily
$\beta_1,\beta'_1 < \beta_2$ and $\gamma_1,\gamma'_1 <
\gamma_3$. Without loss of generality we may assume that $\alpha_1
\leq \alpha'_1$. Then $\lambda_1(u_2^{\beta'_1}u_3^{\gamma'_1}-\lambda_1'
u_1^{\alpha'_1}) -\lambda_1'(u_1^{\alpha'_1-\alpha_1}(u_2^{\beta_1}
u_3^{\gamma_1}-\lambda_1 u_1^{\alpha_1})) = \lambda_1
u_2^{\beta'_1}u_3^{\gamma'_1}-\lambda_1' u_2^{\beta_1}u_3^{\gamma_1}u_1^{
\alpha'_1-\alpha_1 } $.
Factoring out the greatest possible monomial, we obtain an approximate
root of one of the forms (\ref{eqn-approot2}) or (\ref{eqn-approot3}), but
with exponent of $u_3$ strictly less than $\gamma_3$ (respectively,
exponent of $u_2$ strictly less than $\beta_2$), a contradiction.
$\Box$

\begin{rek}
Note that, multiplying each of $u_1,u_2,u_3$ by a suitable non-zero element of
$\R$, we may assume $\lambda_i = 1$ for all $i=1,2,3$.
\end{rek}

Consider a triple of binomials $Q_i,Q_j,Q_k$ with $Q_i =
\sous{u}^{\alpha_i}-\sous{u}^{\beta_i}$, $Q_j =
\sous{u}^{\alpha_j}-\sous{u}^{\beta_j}$, $Q_k =
\sous{u}^{\alpha_k}-\sous{u}^{\beta_k}$,
quasi-homogeneous for a certain $\Q$-weight $\nu_0$, not necessarily approximate
roots.
\medskip

Consider the homomorphism $\sigma : \R[u_1,u_2,u_3] \to \R[[t]]$ defined by
$\sigma(u_q)= t^{\nu_0(u_q)}$, $q=1,2,3$. We have $Q_i,Q_j,Q_k \in \ker(\sigma)$.

\begin{lem}
There exists a syzygy
$\omega_iQ_i+\omega_jQ_j+\omega_kQ_k=0$ where $\omega_i, \omega_j, \omega_k$ are
quasi-homogeneous polynomials in $u_1,u_2,u_3$ with $\omega_i, \omega_j,
\omega_k \notin \ker(\sigma)$.
\end{lem}

\noi Proof :
Let $\nu_0(u_i)=a_i \in \Q$, for $i=1,2,3$. 
Write $\nu_0(\sous{u}^{\alpha_i})=
a_1 \alpha_{i1}+a_2 \alpha_{i2}+a_3 \alpha_{i3}$ and the same for
$\beta_i$, so that $\alpha_i,\beta_i$ belong to the plane
$a_1 x+a_2 y+a_3 z=\nu_0(\sous{u}^{\alpha_i})$.

In the same way, $\alpha_j,\beta_j$ belong to the plane
$a_1 x+a_2 y+a_3 z=\nu_0(\sous{u}^{\alpha_j})$ and $\alpha_k,\beta_k$
belong to the plane $a_1 x+a_2 y+a_3 z=\nu_0(\sous{u}^{\alpha_k})$.

So that the vectors $v_i=\alpha_i-\beta_i$, $v_j=\alpha_j-\beta_j$,
$v_k=\alpha_k-\beta_k$ belong to the plane $a_1 x+a_2 y+a_3 z=0$ in $\Q^3$. So
there is a rational relation of the form $\mu_iv_i+\mu_jv_j+\mu_kv_k=0$.
Multiplying by some integer, we may choose the $\mu_i,\mu_j,\mu_k \in \Z$.

This gives 3 relations between the coordinates :
 \begin{eqnarray}
 \mu_i(\alpha_{i1}-\beta_{i1})+\mu_j(\alpha_{j1}-\beta_{j1})+\mu_k(\alpha_{k1}
-\beta_{k1}) =0  \\
\mu_i(\alpha_{i2}-\beta_{i2})+\mu_j(\alpha_{j2}-\beta_{j2})+\mu_k(\alpha_{k2}
-\beta_{k2}) =0\\
\mu_i(\alpha_{i3}-\beta_{i3})+\mu_j(\alpha_{j3}-\beta_{j3})+\mu_k(\alpha_{k3}
-\beta_{k3}) =0
 \end{eqnarray}

From which we deduce that
\begin{equation}
\left(\frac{(\sous{u}^{\alpha_i})^{\mu_i}}{(\sous{u}^{\beta_i})^{\mu_i}}
\right) \times
\left(\frac{(\sous{u}^{\alpha_j})^{\mu_j}}{(\sous{u}^{\beta_j})^{\mu_j}}\right)
\times \left( \frac{(\sous{u}^{\alpha_k})^{\mu_k}}{(\sous{u}^{\beta_k})^{\mu_k}}
\right) = 1
\end{equation}

and consequently
$$(\sous{u}^{\alpha_i})^{\mu_i}(\sous{u}^{\alpha_j})^{\mu_j}(\sous{u}^{
\alpha_k})^{\mu_k} =
(\sous{u}^{\beta_i})^{\mu_i}(\sous{u}^{\beta_j})^{\mu_j}(\sous{u}^{\beta_k}
)^{\mu_k}$$

Which we can rewrite as
$$
(\sous{u}^{\alpha_i})^{\mu_i}(\sous{u}^{\alpha_j})^{\mu_j}(\sous{u}^{
\alpha_k})^{\mu_k} -
(\sous{u}^{\beta_i})^{\mu_i}(\sous{u}^{\beta_j})^{\mu_j}(\sous{u}^{\beta_k}
)^{\mu_k} = 0.
$$

This last expression can be put under the following form, whatever the sign of
the $\mu_i,\mu_j,\mu_k$ :

\begin{eqnarray*}
\left[(\sous{u}^{\alpha_i})^{\mu_i}-(\sous{u}^{\beta_i})^{\mu_i}\right](\sous{u
}^ { \alpha_j } )^{\mu_j} (\sous{u}^{\alpha_k})^{\mu_k} +
\left[(\sous{u}^{\alpha_j})^{\mu_j}-(\sous{u}^{\beta_j})^{\mu_j}\right]
(\sous{u}^{\beta_i})^{\mu_i}(\sous{u}^{\alpha_k})^{\mu_k} \\  +
\left[(\sous{u}^{\alpha_k})^{\mu_k}-(\sous{u}^{\beta_k})^{\mu_k}\right](\sous{u}
^ { \beta_i })^{\mu_i} (\sous{u}^{\alpha_j})^{\mu_j} =0.
\end{eqnarray*}

Now the relation $a^k-b^k=(a-b)(a^{k-1} + a^{k-2}b + \cdots + b^{k-1})$ applied
to the first bracket shows that
$\left[(\sous{u}^{\alpha_i})^{\mu_i}-(\sous{u}^{\beta_i})^{\mu_i}\right] =
Q_i\times \phi_i$ where $\phi_i$ is a quasi-homogeneous polynomial which is
clearly not in $\ker(\sigma)$. And the same with the two other brackets. This
ends the proof.
$\Box$
\medskip

\noi \textbf{Notation :}  If $Q=\sous{u}^\eta-\lambda \sous{u}^\theta$ is an
approximate root, we denote by $Q'$ the expression
\begin{equation} \label{Qprime}
  Q'=\frac{\sous{u}^\eta}{\sous{u}^\theta}-\lambda.
\end{equation}

Let $G = \oplus_{\gamma \in \Gamma} G_\gamma$ be a graded algebra without zero
divisors. The
\textbf{saturation} of $G$, denoted by $G^*$, is the graded algebra
$$G^* =\{\frac{g}{h} \ |\ g,h \in G, \ h\ homogeneous, \ h\neq 0\}.$$

Assume that $G_\gamma \cong \R$ for all $\gamma \in \Gamma$. Given $\gamma \in
\Gamma$ and $f,g \in G_\gamma,\
g\neq 0$, the notation $\frac{f}{g}$ will mean the unique real number $\lambda$
such that $\lambda g=f$. The real number $\frac{f}{g}$ is independent of the
choice of the isomorphism \begin{equation} \label{eq:congr}
                            G_\gamma \cong \R.
                          \end{equation}
Note that
the number $\lambda$ can be interpreted as an element of $G^*_0 \cong \R$.
\medskip

Now let $\alpha,\beta \in \sper\ A$. Let $<\alpha,\beta>$ be the separating
ideal. Let $\mu_\alpha=\nu_\alpha(<\alpha,\beta>)$ and
$\mu_\beta=\nu_\beta(<\alpha,\beta>)$. Let $Q_i,Q_j$ be two common
approximate roots of $\alpha,\beta$ such that $Q_i,Q_j \notin <\alpha,\beta>$.
Note that, since $\frac{A}{\mathfrak{m}} \cong \R$ and $\alpha,\beta$ are
centered at $\mathfrak{m}$, the graded algebras $\gr_\alpha(A)$ and
$\gr_\beta(A)$ satisfy the condition (\ref{eq:congr}).
\medskip
\begin{deft} \label{def-comp}
We say that the two approximate roots $Q_i,Q_j$ are ($\alpha,\beta$)-comparable
if one of the
following conditions holds :

$$\nu_\alpha(Q'_i)< \nu_\alpha(Q'_j)
\mbox{ and } \nu_\beta(Q'_i)< \nu_\beta(Q'_j)$$

$$\nu_\alpha(Q'_i)>\nu_\alpha(Q'_j)
\mbox{ and } \nu_\beta(Q'_i)> \nu_\beta(Q'_j)$$

$$ \nu_\alpha(Q'_i)= \nu_\alpha(Q'_j)
\mbox{, } \nu_\beta(Q'_i)= \nu_\beta(Q'_j) \mbox{ and }
\frac{\init_\alpha Q'_i}{\init_\alpha Q'_j}
= \frac{\init_\beta Q'_i}{\init_\beta Q'_j}.$$

We say that they are strongly comparable if, up to interchanging $i$
and $j$, we have
\begin{equation}
\nu_\alpha(Q_i)+
\nu_{\alpha 0}(Q_j) < \mu_\alpha\label{eq:stronglycomparable}
\end{equation}
 where $\nu_{\alpha 0}$ is the monomial
valuation such that $\nu_{\alpha 0}(x) = \nu_\alpha(x)$, $\nu_{\alpha 0}(y) =
\nu_\alpha(y)$, $\nu_{\alpha 0}(z) = \nu_\alpha(z)$ (which implies
that (\ref{eq:stronglycomparable}) also holds with $\alpha$ replaced by $\beta$).
\end{deft}

\begin{rek} \label{rek-comp}
Note that ``strongly comparable'' implies ``comparable''. Indeed, write
$Q_i=\sous{u}^{\alpha_i}-\sous{u}^{\beta_i}$,
$Q_j=\sous{u}^{\alpha_j}-\sous{u}^{\beta_j}$. Without loss of generality, we
may write
$Q'_i=\frac{Q_i}{\sous{u}^{\alpha_i}}$, $Q'_j=\frac{Q_j}{\sous{u}^{\alpha_j}}$.
By definition of strongly comparable,
$Q_i\sous{u}^{\alpha_j},Q_j\sous{u}^{\alpha_i} \notin
<\alpha,\beta>$. Thus either
$\nu_\alpha(Q_i\sous{u}^{\alpha_j})<\nu_\alpha(Q_j\sous{u}^{\alpha_i})$,
$\nu_\beta(Q_i\sous{u}^{\alpha_j})<\nu_\beta(Q_j\sous{u}^{\alpha_i})$ or
$\nu_\alpha(Q_i\sous{u}^{\alpha_j})>\nu_\alpha(Q_j\sous{u}^{\alpha_i})$,
$\nu_\beta(Q_i\sous{u}^{\alpha_j})>\nu_\beta(Q_j\sous{u}^{\alpha_i})$ or
$\nu_\alpha(Q_i\sous{u}^{\alpha_j})=\nu_\alpha(Q_j\sous{u}^{\alpha_i})$,
$\nu_\beta(Q_i\sous{u}^{\alpha_j})=\nu_\alpha(Q_j\sous{u}^{\alpha_i})$ and
$\frac{\init_\alpha (Q_i\sous{u}^{\alpha_j})}{\init_\alpha
(Q_j\sous{u}^{\alpha_i})}
= \frac{\init_\beta (Q_i\sous{u}^{\alpha_j})}{\init_\beta
(Q_j\sous{u}^{\alpha_i})}$.
\end{rek}

We saw that there were at most 3 approximate roots of complexity 1 for
$\nu_\alpha$
(and $\nu_\beta$). Suppose there are three such approximate roots common to
$\alpha$ and $\beta$, not in $<\alpha,\beta>$, and denote them by $Q_4,Q_5,Q_6$.
\medskip

%

\begin{lem} $Q_4,Q_5,Q_6$ are either all pairwise comparable or all pairwise
incomparable.
\end{lem}

\noi Proof : Assume that 2 of those roots, say $Q_4$ and $Q_5$, are
comparable. Consider the syzygy
$\omega_4Q_4+\omega_5Q_5+\omega_6Q_6=0$ where $\omega_4,\omega_5,\omega_6$ are
quasi-homogeneous polynomials, not belonging to $\ker(\sigma)$.
Note that this syzygy implies that
\begin{equation} \label{nu0}
 \nu_{\alpha0}(\omega_4 Q_4) = \nu_{\alpha0}(\omega_5 Q_5) =
\nu_{\alpha0}(\omega_6 Q_6)
\end{equation}

We will
prove that $Q_6$ is comparable to $Q_4$ and so, by symmetry, the same
will hold for $Q_6$ and $Q_5$. The following cases are possible :

1. $\nu_\alpha(Q'_4) < \nu_\alpha(Q'_5)$. Then
$\nu_\beta(Q'_4)<\nu_\beta(Q'_5)$. Now using the syzygy and (\ref{nu0}), we
obtain that $\nu_\alpha(\omega_4 Q_4) = \nu_{\alpha 0}(\omega_4 Q_4) +
\nu_\alpha(Q'_4) = \nu_{\alpha 0}(\omega_5 Q_5) + \nu_\alpha(Q'_4) <
\nu_{\alpha0}(\omega_5 Q_5) + \nu_\alpha(Q'_5) = \nu_\alpha(\omega_5 Q_5)$.
So $\nu_\alpha(\omega_4 Q_4) = \nu_\alpha(\omega_6 Q_6)$ which implies that
$\nu_\alpha(Q'_4)=\nu_\alpha(Q'_6)$ and, similarly,
$\nu_\beta(Q'_6)=\nu_\beta(Q'_4)$.

Now, the value of
$\omega_4Q_4+ \omega_6Q_6$ must be greater than the value of each summand
because of the syzygy; in other words, the initial forms of the
summands must cancel each other in the graded algebras of both
$\nu_\alpha$ and $\nu_\beta$.
Hence $\frac{\init_\alpha Q'_6}{\init_\alpha Q'_4}
= \frac{\init_\beta Q'_6}{\init_\beta Q'_4}$, so $Q_6$ is comparable to $Q_4$.

2.  $\nu_\alpha(Q'_4) > \nu_\alpha(Q'_5)$. Then
$\nu_\beta(Q'_4)>\nu_\beta(Q'_5)$, so by symmetry with the previous case
$Q_6$ is comparable to $Q_5$.

3.   $\nu_\alpha(Q'_4) = \nu_\alpha(Q'_5)$,
$\nu_\beta(Q'_4)=\nu_\beta(Q'_5)$ and
\begin{equation}  \label{eq:sauts}
 \frac{\init_\alpha Q'_4}{\init_\alpha
Q'_5}
= \frac{\init_\beta Q'_4}{\init_\beta Q'_5}.
\end{equation}
It follows from $\nu_\alpha(Q'_4) = \nu_\alpha(Q'_5)$ that
$\nu_\alpha(\omega_4Q_4) = \nu_\alpha(\omega_5Q_5)$ and
similarly for $\beta$. Let $\gamma_\alpha=\nu_\alpha(\omega_4Q_4)$ and
$\gamma_\beta=\nu_\beta(\omega_4Q_4)$.

Consider the natural homomorphism of graded
algebras $$\xymatrix{\gr_{\alpha 0}A^* \ar[r]^{\sigma_\alpha}& \gr_\alpha
A^*}.$$
Let $B_\alpha = \sigma_\alpha(\gr_{\alpha 0}A^*)$ and similarly for $\beta$.

Since there are three approximate roots of complexity 1, common to $\alpha$ and
$\beta$, there are at least two $\Q$-linearly independent $\Q$-linear
dependence relations among $\nu_\alpha(x), \nu_\alpha(y),\nu_\alpha(z)$, valid
also for  $\nu_\beta(x), \nu_\beta(y),\nu_\beta(z)$. Thus there is a natural
graded isomorphism
$\xymatrix{gr_{\nu_{\alpha 0}}A \ar[r]^{\iota} & \gr_{\nu_{\beta 0}}A}$. Since
$\iota(\ker \sigma_\alpha)= \ker (\sigma_\beta)$, the map $\iota$ induces an
isomorphism $B_\alpha \to B_\beta$. We obtain the following diagram :
$$\xymatrix{\gr_{\alpha 0}A^* \ar[rr]^{\sigma_\alpha} \ar[rd] \ar[ddd]^\iota &
& \gr_\alpha A^* \\
& B_\alpha \ar@{^{(}->}[ru] \ar[d]^-{\iota_B} \\
& B_\beta \ar@{^{(}->}[rd] & \\
\gr_{\beta 0}A^* \ar[rr]^{\sigma_\alpha} \ar[ru] &  & \gr_\beta
A^*}.$$


If $\omega_i=
M_1^i + \cdots +M_{s_i}^i$ is quasi-homogeneous and not in $\ker(\sigma)$,
we have, for all $j \in \{1,\ldots,s_i\}$,
$\iota(\init_{\alpha 0}(M_j^i))=\init_{\beta 0}(M_j^i)$ and
$\init_{\alpha 0}(\omega_i)=
\init_{\alpha 0}(M_1^i) + \cdots + \init_{\alpha 0}(M^i_{s_i})$ and similarly
for $\beta$. From which we deduce that $\iota(\init_{\alpha 0}(\omega_4))
=\init_{\beta 0}(\omega_4)$. On the other hand, we have
$\iota(\init_{\alpha 0}(\sous{u}^{\alpha_4})) =
\init_{\beta 0}(\sous{u}^{\alpha_4})$ and the same with $\sous{u}^{\alpha_4}$
replaced by $\sous{u}^{\alpha_5}$.
\smallskip

Now we have
$$ \iota\left(
\frac{\init_{\alpha 0}(\sous{u}^{\alpha_4}\omega_4)}{
\init_{\alpha 0}(\sous{u}^{\alpha_5}\omega_5) } \right)=
\frac{\init_{\beta 0}(\sous{u}^{\alpha_4}\omega_4)}{
\init_{\beta 0}(\sous{u}^{\alpha_5}\omega_5) }.$$
Passing to the images in $B_\alpha$ and $B_\beta$ and taking into account that
$\init_{\alpha 0}(\sous{u}^{\alpha_4}\omega_4) \notin \ker(\sigma_\alpha)$, we
obtain the equality of non-zero real numbers
\begin{equation} \label{eq:init}
\frac{\init_{\alpha}(\sous{u}^{\alpha_4}\omega_4)}{
\init_{\alpha}(\sous{u}^{\alpha_5}\omega_5) } =
\frac{\init_{\beta}(\sous{u}^{\alpha_4}\omega_4)}{
\init_{\beta}(\sous{u}^{\alpha_5}\omega_5) }.
\end{equation}

Multiplying the equations (\ref{eq:sauts}) and (\ref{eq:init}), we obtain
\begin{equation}
\frac{\init_\alpha(\omega_4)\init_\alpha(\sous{u}^{\alpha_4})
\init_\alpha(Q'_4)}{\init_\alpha(\omega_5)\init_\alpha(\sous{u}^{\alpha_5})
\init_\alpha(Q'_5)} =
\frac{\init_\beta(\omega_4)\init_\beta(\sous{u}^{\alpha_4})
\init_\beta(Q'_4)}{\init_\beta(\omega_5)\init_\beta(\sous{u}^{\alpha_5})
\init_\beta(Q'_5)}.
\end{equation}
In other words, $$ \frac{\init_\alpha(\omega_4Q_4)}{\init_\alpha(\omega_5Q_5)}
= \frac{\init_\beta(\omega_4Q_4)}{\init_\beta(\omega_5Q_5)}.$$

3a. $\init_\alpha \omega_4Q_4 + \init_\alpha \omega_5Q_5 =0$.
Then
$\init_\beta \omega_4Q_4 + \init_\beta \omega_5Q_5 =0$. Then
$\nu_\alpha(Q'_6) > \nu_\alpha(Q'_4)$ and $\nu_\beta(Q'_6) >
\nu_\beta(Q'_4)$, so $Q_6$ is comparable to $Q_4$.

3b. $\init_\alpha \omega_4Q_4 + \init_\alpha \omega_5Q_5 \neq 0$. Then
$\init_\beta \omega_4Q_4 + \init_\beta \omega_5Q_5 \neq 0$. Then,
using (\ref{nu0}), $\nu_\alpha(Q'_6) = \nu_\alpha(Q'_4)$ and $\nu_\beta(Q'_6) =
\nu_\beta(Q'_4)$. Write the syzygy in the form
$\omega'_4Q'_4+\omega'_5Q'_5+\omega'_6Q'_6=0$ where
$\nu_\alpha(\omega'_4)=\nu_\alpha(\omega'_5)=\nu_\alpha(\omega'_6)$ where, for
$i=4,5,6$, $\omega'_i = \omega_i \sous{u}^{\alpha_i}$.
Now, $$\frac{\init_\alpha Q'_6}{\init_\alpha Q'_4}=\frac{-\init_\alpha
\frac{\omega'_4}{\omega'_6}Q'_4-\init_\alpha
\frac{\omega'_5}{\omega'_6}Q'_5}{\init_\alpha Q'_4} = -\init_\alpha
\frac{\omega'_4}{\omega'_6}- \init_\alpha
\frac{\omega'_5}{\omega'_6}\frac{\init_\alpha Q'_5}{\init_\alpha
  Q'_4} = -\init_\beta \frac{\omega'_4}{\omega'_6}- \init_\beta
\frac{\omega'_5}{\omega'_6}\frac{\init_\beta Q'_5}{\init_\beta
  Q'_4} = \frac{\init_\beta Q'_6}{\init_\beta Q'_4}$$ and  so again
$Q_6$ is comparable to $Q_4$. 
\medskip

\noi \textbf{From now on, assume that $A$ is excellent.}

\subsection{Some of $u_1,u_2,u_3$ belong to the separating ideal}

Assume that $\{u_1,u_2,u_3\} \cap <\alpha,\beta> \neq
\emptyset$. The case
$u_1,u_2,u_3 \in <\alpha,\beta>$ is trivial since then $<\alpha
,\beta> = \mathfrak{m}$. 

If $u_1 \notin <\alpha,\beta>$ and $u_2,u_3 \in <\alpha,\beta>$, then the only
relevant approximate roots appearing in (\ref{eqn-gi}) and
(\ref{eq:dominant}) are $u_1,u_2,u_3$. Then all the inequalities defining $C'$
are monomial in $u_1,u_2,u_3$. So $C'$ is connected by Corollary
\ref{prop-connexe}. 

Finally, suppose that $u_1,u_2 \notin <\alpha, \beta>$,
$u_3 \in <\alpha, \beta>$. Then all the approximate roots belong to $\{u_3\}
\cup \R[u_1,u_2]$. 

After a suitable sequence of affine monomial blowings up  $A \to A'$
such that $\sper\ A'$ contains the common center $\mathfrak{m}'$ of $\alpha'$
and $\beta'$, $\mathfrak{m}'$ has a regular system of parameters
$u_1',u_2',u_3'$ such that all the approximate roots are monomials in
$u_1',u_2',u_3'$ up to multiplication by units of $A'$. Let $\xi'$ be the
unique point of $\sper\ A'$ with support $\mathfrak{m}'$. 

Let $\pi: \sper\ A' \to \sper\ A$ the induced map of real spectra. Write
$\mathbf{Q}^{\delta_i} = \sous{u}'^{\delta'_i} v_i$, $\mathbf{Q}^{\delta_{ji}}
= \sous{u}'^{\delta'_{ji}} v_{ji}$ with $v_i, v_{ji} \in A'\setminus
\mathfrak{m}'$. 
If $C$ is as in (\ref{ensCC}), the set $\pi^{-1}(C)$ contains the
set $$\tilde{C} = \left\lbrace 
\begin{array}{l} 
\delta \text{ is centered at } \mathfrak{m}' \\ 
\nu_\delta(\sous{u}'^{\delta'_i}) < \nu_\delta(\sous{u}'^{\delta'_{ji}}),
\ i \in \{1,\ldots,s\}, j \in \{1,\dots, N_i\} \\ 
sgn_\delta u'_\ell = sgn_\alpha u'_\ell \text{ for all } u'_\ell \text{
appearing in } \sous{u}'^{\delta'_i} \text{ for some } i 
\end{array}
\right\rbrace.$$
Take $d \in \R$ such that $\dsp{d > \max_{1\leq i \leq s} N_i \times
\max_{\begin{array}{c}1 \leq i \leq s\\1 \leq j \leq
N_i \end{array}}\frac{|v_{ji}(\xi')|}{|v_i(\xi')|}}$.

Let $$\tilde{\tilde{C}} = \left\lbrace 
\begin{array}{l} 
\delta \text{ is centered at } \mathfrak{m}' \\ 
|\sous{u}'^{\delta'_i}(\delta)| > d |\sous{u}'^{\delta'_{ji}}(\delta)| 
\ i \in \{1,\ldots,s\}, j \in \{1,\dots, N_i\} \\ 
sgn_\delta u'_\ell = sgn_\alpha u'_\ell \text{ for all } u'_\ell \text{
appearing
in } \sous{u}'^{\delta'_i} \text{ for some } i 
\end{array}
\right\rbrace.$$

$\tilde{\tilde{C}}$ is connected by Corollary \ref{prop-connexe}. 

Then $\alpha',\beta' \in \tilde{\tilde{C}}$, hence $\alpha,\beta \in
\pi(\tilde{\tilde{C}})$, $\pi(\tilde{\tilde{C}})$ is connected and contained in
$\{g_1\cdots g_s \neq 0\}$. This completes the proof in the case when
$\{u_1,u_2,u_3\} \cap <\alpha,\beta> \neq   \emptyset$. 
\medskip

\textbf{From now on, we assume $\{u_1,u_2,u_3\} \cap <\alpha,\beta> =
\emptyset$ and that, unless otherwise specified, there are 3 relevant
approximate roots $Q_4,Q_5,Q_6$.}

\subsection{All the approximate roots are pairwise comparable} \label{sec:cas1}

Without loss of generality, assume that $\nu_\alpha(Q'_6) \leq
\nu_\alpha(Q'_4)$, $\nu_\alpha(Q'_6) \leq \nu_\alpha(Q'_5)$ and
$u_1(\alpha)>0,u_2(\alpha)>0,u_3(\alpha)>0, Q_6(\alpha)>0$.

Write \begin{equation}\label{eqn-qu6}
Q_6= -\frac{\omega_4}{\omega_6}Q_4 -
\frac{\omega_5}{\omega_6}Q_5.
\end{equation}

Assume that
$\frac{\omega_5}{\omega_6}Q_5>_{\alpha,\beta}0 $
(if the opposite holds, a similar
reasoning applies). Then $Q_6 >_{\alpha,\beta}0 \Ra
|\frac{\omega_4}{\omega_6}Q_4|
>_{\alpha,\beta} |\frac{\omega_5}{\omega_6}Q_5|$. There exists
$\epsilon >0$ in $\R$ such that
\begin{equation} \label{eqn-ineq1}
(1-\epsilon)|\frac{\omega_4}{\omega_6}Q_4|
>_{\alpha,\beta} |\frac{\omega_5}{\omega_6}Q_5|.
\end{equation}
Then
\begin{equation}
|Q_6|>_{\alpha,\beta} \epsilon |\frac{\omega_4}{\omega_6}Q_4|. \label{eqn-ineq2}
\end{equation}
\medskip

We now describe the connected set required in the Connectedness
Conjecture which we will define by inequalities among
certain generalized monomials and sign conditions on the $Q_i$. Let
$g_1,\ldots,g_s \in A$ be as in the statement of the conjecture. Let
\begin{equation} 
g_i = \mathbf{Q}^{\delta_i} +
\sum_{j=1}^{N_i}c_{ji}\mathbf{Q}^{\delta_{ji}},\ i\in \{1,\ldots,s\}
\end{equation}
be the standard form of $g_i$.
\medskip

For each $i \in \{1,\ldots,s\}$, in the sum
$\sum_{j=1}^{N_i}c_{ji}\mathbf{Q}^{\delta_{ji}}$ replace $Q_6$ by the
right hand side of (\ref{eqn-qu6}) and write the result as a sum of
generalized monomials (with possibly negative exponents) in
$u_1,u_2,u_3,Q_4,Q_5$ : $$\sum_{j=1}^{N'_i} c'_{ji}\mathbf{Q}^{\delta'_{ji}}.$$

In each generalized monomial $\mathbf{Q}^{\delta_i}$, replace $Q_6$
by $\epsilon \frac{\omega_4}{\omega_6}Q_4$ and let
$c'_i\mathbf{Q}^{\delta'_i}$ be the resulting generalized monomial.

Let $D$ be the subset of $\sper\ A$ consisting of points $\delta$ such
that 
\begin{eqnarray}
|c'_i\mathbf{Q}^{\delta'_i}(\delta)| >
\frac{1}{N'_i}|c'_{ji}\mathbf{Q}^{\delta'_{ji}}(\delta)|& i \in
\{1,\ldots,s\},\ j \in \{1,\ldots,N_i\} \label{eqn-array1}\\
 (1-\epsilon)|\frac{\omega_4}{\omega_6}Q_4(\delta)|
>|\frac{\omega_5}{\omega_6}Q_5(\delta)| & \label{eqn-array2}\\
sgn(u_1(\delta)) =sgn(u_1(\alpha))\label{eqn-array3}& \\ sgn(u_2(\delta))
=sgn(u_2(\alpha))&\\
sgn(u_3(\delta)) =sgn(u_3(\alpha))&\\ sgn(Q_4(\delta)) =sgn(Q_4(\alpha))&\\
sgn(Q_5(\delta)) =sgn(Q_5(\alpha))\label{eqn-array7}&
\end{eqnarray}

By definition of standard form,
$\nu_\alpha(\mathbf{Q}^{\delta_i}) <
\nu_\alpha(\mathbf{Q}^{\delta_{ji}})$ for all $i,j$ and similarly for
$\nu_\beta$. By (\ref{eqn-ineq2}), this implies that
$\nu_\alpha(\mathbf{Q}^{\delta'_i})
<\nu_\alpha(\mathbf{Q}^{\delta'_{ji}})$ for all $i,j$ and the same for
$\nu_\beta$. Thus inequalities (\ref{eqn-array1}) hold for
$\delta=\alpha$ and $\delta=\beta$. This proves that $\alpha,\ \beta
\in D$.
The polynomials $g_i$ have constant sign on $D$ for all $i$ because
the inequalities ensure that the sign of $g_i$ is determined by the
sign of its dominant monomial $\mathbf{Q}^{\delta_i}$. With these
conditions, using (\ref{eqn-qu6}), we see that the signs of both $Q_6$
and $\mathbf{Q}^{\delta_i}$ are constant on $D$.
\medskip

It remains to prove that $\alpha$ and $\beta$ belong to the same connected component of $D$ Let $A \to A'$ be a finite
sequence of affine monomial blowings up such that
$\sper\ A'$ contains the common center $m'$ of $\alpha'$ and
$\beta'$ and there is a regular system of parameters $(x',y',z')$ at
$m'$ such that $x'$ is a monomial in $u_1,u_2,u_3$, $y'=Q'_4$ and $z'=Q'_5$.

For each inequality
\begin{equation} \label{eqn-ineg}
|c\mathbf{Q}^\epsilon(\delta) | < |d\mathbf{Q}^\gamma(\delta)|
\end{equation} appearing in the
definition of $D$, there exist $\epsilon'_x,\epsilon'_y,\epsilon'_z,
\gamma'_x , \gamma'_y,  \gamma'_z \in \Z$  and
elements $u,v \in{A'_{m'}}\setminus m'A'_{m'}$
such that $c\mathbf{Q}^\epsilon =
ux'^{\epsilon'_x}y'^{\epsilon'_y}z'^{\epsilon'_z}$ and $d\mathbf{Q}^\gamma =
vx'^{\gamma'_x}y'^{\gamma'_y}z'^{\gamma'_z}$. Take positive constants
$\tilde{u},\tilde{v} \in R$ such that $|u|<\tilde{u}$ and $|v| >
\tilde{v}$. Then, for any $\delta' \in \sper\
A'_{m'} 
$, the inequality
\begin{equation}\label{eqn-ineg2}
|\tilde{u}x'^{\epsilon'_x}(\delta')y'^{\epsilon'_y}(\delta')z'^{\epsilon'_z}(\delta')|
< |\tilde{v}x'^{\gamma'_x}(\delta')y'^{\gamma'_y}(\delta')z'^{\gamma'_z}(\delta')|
\end{equation} implies the
inequality (\ref{eqn-ineg}) where $\delta$ is the image of $\delta'$
in $\sper\ A$.

Since \begin{equation}\label{eqn-ineg3}
\nu_\alpha(x'^{\epsilon'_x}y'^{\epsilon'_y}z'^{\epsilon'_z}) >
\nu_\alpha(x'^{\gamma'_x}y'^{\gamma'_y}z'^{\gamma'_z})
\end{equation} and similarly for $\beta$, the inequalities (\ref{eqn-ineg2})
hold for both $\delta'=\alpha'$ and $\delta'=\beta'$.
\medskip

Let $D'$ be the subset of $\sper\ A'_{m'}$ 
defined by all the resulting inequalities of the form 
(\ref{eqn-ineg2}) and the sign conditions
\begin{eqnarray}
sgn(x'(\delta')) =sgn(x'(\alpha')) \\ sgn(y'(\delta')) =sgn(y'(\alpha'))\\
sgn(z'(\delta')) =sgn(z'(\alpha')).
\end{eqnarray}

The set $D'$ is connected by Corollary \ref{prop-connexe}. Its image
in $\sper\ A$ is connected, contains $\alpha$ and $\beta$ and is
contained in the set $\{g_1\cdots g_s \neq 0\}$. This completes the
proof of the Connectedness Conjecture in the case when
$Q_4,Q_5,Q_6$
are pairwise comparable.
\bigskip

\begin{rek} \label{rek_limit}
The same method works to prove the Connectedness conjecture also in the case
when $\nu_\beta(Q'_4)=\nu_\beta(Q'_5)=\nu_\beta(Q'_6)$ and
$\nu_\alpha(Q'_4)=\nu_\alpha(Q'_5)<\nu_\alpha(Q'_6)$.
\end{rek}


\subsection{The case when there are only one or two relevant approximate
roots}

Assume that there are exactly 2 approximate roots $Q_4,Q_5\notin<\alpha,\beta>$, common to $\alpha$ and $\beta$. We proceed as in the case of 3 comparable approximate roots. This means that,
after a suitable sequence of affine monomials blowings up  $A \to A'$
such that $\sper\ A'$ contains the common center $m'$ of $\alpha'$ and
$\beta'$, $m'$ has a regular system of parameters $x',y',z'$ such that $x'$ is
a monomial in $u_1,u_2,u_3$, $y'=Q'_4$, $z'=Q'_5$ (remember the
notation (\ref{Qprime})).

Then as above, we replace the inequalities of type (\ref{eqn-ineg}) by
inequalities involving only monomials as in (\ref{eqn-ineg2}). Once
again we apply the Corollary \ref{prop-connexe} to ensure the existence of a
set $\mathcal{C}$ as required.

The case with only one relevant approximate root, $Q_4$, is more difficult.
\medskip

\noi \textbf{Claim.} There exists a connected subset $\tilde{C}$ of $C$, which we
will describe explicitly, containing $\alpha$ and $\beta$.
\medskip

\noi Proof : First, we consider the special case when $A$ is the localization
of the polynomial ring $A=\R[x,y,z]_{(x,y,z)}$. Let $A \to A_1 \to \cdots \to
A_i$ be a finite sequence of affine monomial blowings up with respect to
$\alpha$
(see \cite{LMSS}, Proposition 6.1) such that $\sous{u}^{\alpha_4-\beta_4} \in
A_i$ where $\sous{u}=(x,y,z)$. Let $z_i =
\sous{u}^{\alpha_4-\beta_4}-1$. By construction, $A_i$ is the
localization of a polynomial ring of the form $\R[x_i,y_i,z_i]$, where $x_i,y_i$
are monomials in $x,y,z$ and $z_i=\sous{u}^{\alpha_4-\beta_4}-1$, by the
multiplicative set $\R[x,y,z]\setminus (x,y,z)$. Let $\sous{u}_i=(x_i,y_i,z_i)$.

In $A_i$ the inequalities (\ref{eq:dominant}) can be rewritten as
\begin{equation}
 \left| (c_i+z_if_i)\sous{u}_i^{\gamma_i} \right| < \left|
(d_i+z_ih_i)\sous{u}_i^{\delta_i} \right| \text{ where } c_i,d_i \in \R,\ c_id_i
\neq 0,\ f_i,h_i \in \R[x_i,y_i,z_i].
\end{equation}

Let $\tilde{C} \subset \sper\ A_i$ be the set defined by the stronger
inequalities
\begin{equation}\label{eq:Ctilde}
 \left| 2c_i\sous{u}_i^{\gamma_i} \right| \leq \left|
\frac{1}{2}d_i\sous{u}_i^{\delta_i} \right|
\end{equation}
and the sign conditions $sgn(x_i(\delta))= sgn(x_i(\alpha))$, $sgn(y_i(\delta))
= sgn(y_i(\alpha))$, $sgn(z_i(\delta)) = sgn(z_i(\alpha))$. The set $\tilde{C}$
contains the strict transforms of $\alpha$ and $\beta$.

Let $D$ be the subset of $\sper\ \R[x_i,y_i,z_i]$ defined by the inequalities
(\ref{eq:Ctilde}) and the same sign conditions as above.

Using the cartesian diagram
$$\xymatrix{\sper\ \R[x,y,z] & & \ar[ll] \sper\ \R[x_i,y_i,z_i]  \\ & \Box & \\
\sper\ \R[x,y,z]_{(x,y,z)} \ar@{^{(}->}[uu] & & \sper\ \R[x_i,y_i,z_i]_S
\ar@{^{(}->}[uu] \ar[ll] }$$ where $S=\R[x,y,z]\setminus (x,y,z)$,
we can identify
$\tilde{C}$ with $$\bigcap_{N=1}^\infty D \cap \{\delta \in \sper\ \R[x_i,y_i,
z_i] \ ;\ |x(\delta)| \leq \frac{1}{N},\ |y(\delta)| \leq \frac{1}{N},\
|z(\delta)| \leq \frac{1}{N} \}.$$ By Lemma 4.1 of \cite{LMSS}, we have
that $D \cap \{\delta \in \sper\ \R[x_i,y_i,
z_i] \ ;\ |x(\delta)| \leq \frac{1}{N},\ |y(\delta)| \leq \frac{1}{N},\
|z(\delta)| \leq \frac{1}{N} \}$ is connected, so applying Lemma 7.1 of
\cite{LMSS}, we deduce that the intersection is a non empty closed
connected set, so $\tilde{C}$ is connected and contains the strict transforms
of $\alpha$ and $\beta$, hence its image in $\sper\ A$ is the desired connected
set. This completes the proof when $A=\R[x,y,z]_{(x,y,z)}$. The general case now
follows from Theorem
\ref{th:connex}.

\subsection{The approximate roots are pairwise incomparable}
\begin{lem}
At least two of $Q_4,Q_5,Q_6$ have $\nu_\alpha$-value strictly greater
than $\frac{\mu_\alpha}{2}$ (and similarly for $\nu_\beta$-value).
\end{lem}

\noi Proof : This is an immediate consequence of Definition
\ref{def-comp} and Remark \ref{rek-comp}. $\Box$
\medskip

Without loss of generality, assume that
\begin{equation} \label{eqn-order}
\nu_\alpha(Q_4) \leq
\nu_\alpha(Q_5) \leq \nu_\alpha(Q_6). \end{equation}
Then
\begin{equation} \label{eqn-convention}
\nu_\alpha(Q_5) > \frac{\mu_\alpha}{2},\ \nu_\alpha(Q_6) >
\frac{\mu_\alpha}{2}.
\end{equation}
\begin{prop}
Consider a generalized monomial $\mathbf{Q}^\gamma$ divisible by one
of $Q_4Q_5$, $Q_4Q_6$, $Q_5Q_6$, $Q_5^2$, $Q_6^2 $.

\noi Then (a) $ \mathbf{Q}^\gamma \in <\alpha,\beta>$

(b) $\mathbf{Q}^\gamma$ belongs to the ideal generated by all the
generalized monomials belonging to\linebreak $<\alpha,\beta>$ and not divisible
by any of $Q_5^2,Q_6^2,Q_4Q_5,Q_4Q_6,Q_5Q_6$.
\end{prop}

\noi Proof : (a) If $\left.Q_5^2\ \right|\ \mathbf{Q}^\gamma$ or
$\left.Q_6^2\ \right|\ \mathbf{Q}^\gamma$, the result follows immediately from
(\ref{eqn-convention}).

\noi If $Q_4Q_5\ \left|\ \mathbf{Q}^\gamma\right.$, $Q_4Q_6\ \left|\ \mathbf{Q}^\gamma\right.$ or
$Q_5Q_6\ \left|\ \mathbf{Q}^\gamma\right.$, the result follows from Definition
\ref{def-comp} and Lemma \ref{rek-comp}.

(b) By pairwise incomparability and (\ref{eqn-order}), we have
$\nu_\alpha(Q_5)+\nu_{0\alpha}(Q_5) \geq \mu_\alpha$ and similarly for
$\beta$. As well, \begin{equation} \label{eqn-loc}
\nu_\alpha(Q_6)+\nu_{0\alpha}(Q_6) \geq
\mu_\alpha.\end{equation} Suppose, for example, $\left.Q_5^2\ \right|\ \mathbf{Q}^\gamma$. Write
$Q_5 = \omega-\epsilon$ where $\omega,\epsilon$ are monomials in
$x,y,z$.
Then $\mathbf{Q}^\gamma$ belongs to the ideal generated by
$Q_5\epsilon$, $Q_5\omega$ and by (\ref{eqn-loc}), $Q_5\epsilon$,
$Q_5\omega \in <\alpha,\beta>$. The cases when $\mathbf{Q}^\gamma$ is
divisible by $Q_6^2$, $Q_4Q_5$, $Q_4Q_6$, $Q_5Q_6$ are handled
similarly. $\Box$
\medskip

We now describe the connected set required in the Connectedness
Conjecture by inequalities on the size of
certain generalized monomials and sign conditions on the $Q_i$. Let
$g_1,\ldots,g_s \in A$ be as in the statement of the conjecture. Let
\begin{equation}\label{eqn-gi2}
g_i = \mathbf{Q}^{\delta_i} +
\sum_{j=1}^{N_i}c_{ji}\mathbf{Q}^{\delta_{ji}},\ i\in \{1,\ldots,s\}
\end{equation} be the standard form of $g_i$ of level
$\nu_\alpha(<\alpha,\beta>$.
\medskip

Let $\mathcal{S}$ be a finite set of generalized monomials, not divisible by
$Q_5^2,Q_6^2$, $Q_4Q_5$, $Q_4Q_6$, $Q_5Q_6$, belonging
to $<\alpha,\beta>$, which generate $<\alpha,\beta>$. In addition, we require all the
monomials $\mathbf{Q}^\lambda \in \mathcal{S}$ to have the following
property : if $Q_4\ \left|\ \mathbf{Q}^\lambda\right.$ then
\begin{equation}
\nu_\alpha(\mathbf{Q}^\lambda) - \nu_\alpha(Q_4) + \nu_{\alpha 0}(Q_4)
< \mu_\alpha
\end{equation}
and similarly for $Q_5$ and $Q_6$.

Let
$\mathcal{T}$ be the set of all the generalized monomials not belonging to
$<\alpha,\beta>$.

Let $C$ be the subset of $\sper\ A$ consisting of points $\delta$ such
that
\begin{eqnarray}
\nu_{\delta}(\mathbf{Q}^\gamma) < \nu_{\delta}(\mathbf{Q}^\lambda) & \label{eqn-array12}\\
\nu_{\delta}(\mathbf{Q}^\theta)=\nu_{\delta}(\mathbf{Q}^\eta) & \label{eqn-array22}\\
sgn(u_1(\delta)) =sgn(u_1(\alpha))\label{eqn-array32}& \\ sgn(u_2(\delta))
=sgn(u_2(\alpha))&\\
sgn(u_3(\delta)) =sgn(u_3(\alpha))&\\ sgn(Q_4(\delta)) =sgn(Q_4(\alpha))&\\
sgn(Q_5(\delta)) =sgn(Q_5(\alpha))\label{eqn-array72}&\\
\sgn(Q_6(\delta)) =sgn(Q_6(\alpha)) \label{eqn-array82} &
\end{eqnarray}
where
(i) $\mathbf{Q}^\theta,\mathbf{Q}^\eta$ run over all the pairs of
elements of $\mathcal{T}$ satisfying (\ref{eqn-array22}) for
$\delta=\alpha$ and $\delta=\beta$,

(ii) $\mathbf{Q}^\gamma,\mathbf{Q}^\lambda$ run over all the pairs of
generalized monomials such that $\mathbf{Q}^\gamma \in \mathcal{T}$,
$\mathbf{Q}^\lambda \in \mathcal{T}\cup \mathcal{S}$ and
(\ref{eqn-array12}) holds for
$\delta=\alpha$ and $\delta=\beta$.
\medskip

Note that the definition of $C$ implies that for all $\delta\in C$, the
binomials $Q_4,Q_5,Q_6$ are approximate roots for the valuation $\nu_\delta$.
\medskip

All points $\delta \in C$ share, by definition of $C$, the same
approximate roots $Q_4,Q_5,Q_6$. This implies that the
$\dim_\Q(\sum_{j=1}^3 \Q \nu_\delta(u_i)) =1$, for all $\delta
\in C$ (in particular, for $\delta=\alpha$ and $\delta=\beta$). Moreover, there
exist $r,s \in \Q$ such that $\nu_\delta(u_2)= r\nu_\delta(u_1)$ and
$\nu_\delta(u_3)=s\nu_\delta(u_1)$ for all $\delta
\in C$. Then each
equality or inequality of (\ref{eqn-array12}), (\ref{eqn-array22}) may be
written in a form containing only $\nu(u_1)$ and $\nu(Q_4),\nu(Q_5),\nu(Q_6)$.
As $\nu_{\delta0}(\mathbf{Q}^\lambda)$ can be written purely in terms of
$\nu(u_1)$ and
as $\nu(Q_\ell) = \nu_0(Q_\ell)+\nu(Q'_\ell)$, any relation of the
form (\ref{eqn-array12}) or (\ref{eqn-array22}) may be written in terms
of $\nu(u_1),\nu(Q_4'),\nu(Q_5'),\nu(Q_6')$.

\begin{prop} \label{prop-diag}
$C$ contains a point $\epsilon$ such that
$\nu_\epsilon(Q'_4)=\nu_\epsilon(Q'_5)=\nu_\epsilon(Q'_6)$.
\end{prop}

\noi Without loss of generality, assume that
\begin{eqnarray}
\nu_\alpha(Q'_4) =\nu_\alpha(Q'_5) < \nu_\alpha(Q'_6)
&\mbox{ and } \\
 \nu_\beta(Q'_4) =\nu_\beta(Q'_6) < \nu_\alpha(Q'_5)&
\end{eqnarray}
  Replacing $\alpha$ by $\alpha'$ lying in $C$ such that
\begin{eqnarray}
\nu_{\alpha'}(Q'_4) =\nu_{\alpha'}(Q'_5) < \nu_{\alpha'}(Q'_6)
\end{eqnarray} and $\Gamma_{\alpha'} \subset \R$ does not change the
problem and similarly for $\beta$. From now on, we will assume that
$\Gamma_\alpha \subset \R$ and $\Gamma_\beta\subset \R$.
\medskip

Let 
\begin{equation}
 \phi : \left\{ \delta \in \sper\ A\ \left| 
\begin{array}{l}\ \Gamma_\delta \subset \R \\
\delta \text{ centered in } \mathfrak{m} \\ 
Q_4,Q_5,Q_6 \text{ are approximate } \\ \text{roots for } \delta . 
\end{array} \right.\right\} \to \begin{array}{l}\{(a_1,a_2,a_3,a_4) \in \R^4\\
\ \ \ |\ a_1,a_2,a_3,a_4 >0 \} \end{array}
\end{equation} 
be the map defined by $\phi(\delta) = 
(\nu_\delta(u_1),\nu_\delta(Q'_4),\nu_\delta(Q'_5),\nu_\delta(Q'_6))$. 

\begin{lem} \label{lem:iff}
 A point $(a_1,a_2,a_3,a_4) \in \R^4$, $a_1>0,a_2>0,a_3>0,a_4>0$ is in the
image of $\phi$ if and only if one of the following conditions holds 
\begin{equation} \label{conditions} a_2 = a_3
\leq a_4, \ a_2 = a_4 \leq a_3, \ a_3 = a_4 \leq a_2.
\end{equation}
\end{lem}

\noi Proof : Write the syzygy in the form
$\omega'_4Q'_4+\omega'_5Q'_5+\omega'_6Q'_6=0$. The ``only if'' part follows
from this.

``If'' : Suppose, for example, that $a_2=a_3 \leq a_4$. Consider a sequence of
blowings up $A \to A'$ such that $A'$ has a maximal ideal $m'$ with a regular
system of parameters $(x',y',z')$ such that $Q'_4=y',Q'_6=z'$ and 
$u_1,u_2,u_3$ are monomials in $x'$, up to multiplication by units of $A'$.
Write $ u_1=x'^\gamma v$ with $\gamma \in \N^*_+$, $v \in A'\setminus m'$. Take
a point $\delta \sper\ A$ such that $\nu_\delta(x')=\frac{a_1}{\gamma}$,
$\nu_\delta(y')=a_2$, $\nu_\delta(z')=a_4$. Then $\delta$ is centered in
$\mathfrak{m}$ and has $Q_4,Q_5,Q_6$ as approximate roots. this proves that
$(a_1,a_2,a_3,a_4) \in Im(\phi)$. 
\medskip

Next, we reduce the problem to the case when each of the inequalities (\ref{eqn-array12})
and equalities (\ref{eqn-array22}) involves at
most one of $Q_4,Q_5,Q_6$ (possibly raised to some power). Namely let
$h_1,\ldots,h_p$ be the complete list of inequalities (\ref{eqn-array12}) and equalities
(\ref{eqn-array22}). We order the $h_i$ in such a
way that all the inequalities-equalities involving at most one of
$Q_4,Q_5,Q_6$ come first and those involving at least 2 of
$Q_4,Q_5,Q_6$ come later. In each of the above two lists, we order the
inequalities-equalities by the value of the left-hand side.

\begin{lem}
Assume Proposition \ref{prop-diag} is true in the special case when each $h_i$
contains at most one of $Q_4$, $Q_5$, $Q_6$. Then it is true in general.
\end{lem}

\noi Proof : We argue by contradiction. Suppose that Proposition
\ref{prop-diag} is false for $h_1,\ldots,h_p$.
Without loss of generality, we may assume that Proposition is true for
$h_1,\ldots,h_{p-1}$. Let
$\tilde{C} \supset C$ be the set defined by the same conditions as $C$
except for $h_p$. By assumptions, $\tilde{C}$ contains a point
$\epsilon$ such that
$\nu_\epsilon(Q_4')=\nu_\epsilon(Q_5')=\nu_\epsilon(Q_6')$.
\medskip

In the following formulae, the notation 
\begin{equation} \label{eq:ensemble}
 \{\nu_\delta(\mathbf{Q}^\gamma) = \nu_\delta(\mathbf{Q}^\lambda)\}
\end{equation}
means ``the set of all the points of $\R^4$ of the form
$(\nu_\delta(u_1),\nu_\delta(Q'_4),\nu_\delta(Q'_5),
\nu_\delta(Q'_6))$ where $\delta$ satisfies (\ref{eq:ensemble})''. This
notation makes sense because because $\nu_\delta(\mathbf{Q}^\gamma)$ and
$\nu_\delta(\mathbf{Q}\lambda)$ are completely determined by
$(\nu_\delta(u_1),\nu_\delta(Q'_4),\nu_\delta(Q'_5),\nu_\delta(Q'_6))$.
The set (\ref{eq:ensemble}) is contained in a hyperplane $H$ of $\R^4$, defined
by a linear equation with rational coefficients, and contains the subset of $H$
satisfying the conditions (\ref{conditions}).

If $h_p$ is a strict inequality, write $h_p$ in the form
$\nu_\delta(\mathbf{Q}^\gamma)<\nu_\delta(\mathbf{Q}^\lambda)$ and consider two
segments in $\R^4$
\begin{equation}
[(\nu_\alpha(u_1),\nu_\alpha(Q'_4),\nu_\alpha(Q'_5),\nu_\alpha(Q'_6)),
(\nu_\epsilon(u_1),\nu_\epsilon(Q'_4),\nu_\epsilon(Q'_5),\nu_\epsilon(Q'_6))]
\end{equation}
\begin{equation}
[(\nu_\beta(u_1),\nu_\beta(Q'_4),\nu_\beta(Q'_5),\nu_\beta(Q'_6)),
(\nu_\epsilon(u_1),\nu_\epsilon(Q'_4),\nu_\epsilon(Q'_5),\nu_\epsilon(Q'_6))].
\end{equation}

Since $\nu_\epsilon(Q'_4)=\nu_\epsilon(Q'_5)=\nu_\epsilon(Q'_6)$ and the
left endpoint
of $$[(\nu_\alpha(u_1),\nu_\alpha(Q'_4),\nu_\alpha(Q'_5),\nu_\alpha(Q'_6)),
(\nu_\epsilon(u_1),\nu_\epsilon(Q'_4),\nu_\epsilon(Q'_5),
\nu_\epsilon(Q'_6))] $$
satisfies the conditions (\ref{conditions}), so does every point of that
interval. The same holds for the interval
$$[(\nu_\beta(u_1),\nu_\beta(Q'_4),\nu_\beta(Q'_5),\nu_\beta(Q'_6)),
(\nu_\epsilon(u_1),\nu_\epsilon(Q'_4),\nu_\epsilon(Q'_5),
\nu_\epsilon(Q'_6))].$$

Since $\epsilon \notin C$, the following intersections are non empty; each of
them consists of one point 
\begin{eqnarray}
[(\nu_\alpha(u_1),\nu_\alpha(Q'_4),\nu_\alpha(Q'_5),\nu_\alpha(Q'_6)),
(\nu_\epsilon(u_1),\nu_\epsilon(Q'_4),\nu_\epsilon(Q'_5),
\nu_\epsilon(Q'_6))] \\
\cap \left\{ \nu_\delta(\mathbf{Q}^\gamma) 
=\nu_\delta(\mathbf{Q}^\lambda) \right\}  =: \{(a_1,a_2,a_3,a_4)\}
\end{eqnarray}
and
\begin{eqnarray}
[(\nu_\beta(u_1),\nu_\beta(Q'_4),\nu_\beta(Q'_5),\nu_\beta(Q'_6)),
(\nu_\epsilon(u_1),\nu_\epsilon(Q'_4),\nu_\epsilon(Q'_5),
\nu_\epsilon(Q'_6))] \\
\cap \{\nu_\delta(\mathbf{Q}^\gamma) =\nu_\delta(\mathbf{Q}^\lambda)\}
=: \{(b_1,b_2,b_3,b_4)\}.
\end{eqnarray}
Take points $\alpha_0 \in \phi^{-1}((a_1,a_2,a_3,a_4))$ and $\beta_0 \in
\phi^{-1}((b_1,b_2,b_3,b_4))$.
\medskip

\noi In particular
\begin{eqnarray} \label{eqn:first} \nu_{\alpha_0}(Q'_4) =\nu_{\alpha_0}(Q'_5)
\leq
  \nu_{\alpha_0}(Q'_6) & \text{and} \label{eqn-compare1}\\ \label{eqn:second}
\nu_{\beta_0}(Q'_4) =\nu_{\beta_0}(Q'_6) \leq
  \nu_{\beta_0}(Q'_5)\label{eqn-compare2}
\end{eqnarray}
and at least one of the two inequalities is strict. Note that $\alpha_0 \neq
\beta_0$.
Suppose the inequality (\ref{eqn:first}) is strict. Let $h_{p_0}$ be the
equality $\nu_\delta(\mathbf{Q}^\gamma)=\nu_\delta(\mathbf{Q}^\lambda)$.
\medskip

If $h_p$ is an equality, put $\alpha_0=\alpha$,
$\beta_0=\beta$ and let $h_{p_0}=h_p$. 

We can now contradict (\ref{eqn-compare1}), (\ref{eqn-compare2}) as follows.
Suppose, for example, that $h_{p_0}$ has the form
\begin{equation}
\nu_\delta(Q_4^a \omega)=\nu_\delta(Q_5 \eta)
\end{equation}
$a \geq 1$ and write $Q_4=\epsilon_4-\omega_4$,
$Q_5=\epsilon_5-\omega_5$.

Then $Q_4^{a-1}\omega \epsilon_4$ and $\epsilon_5\eta$ do not
belong to $<\alpha,\beta>$ and the relation between
$\nu_\alpha(Q_4^{a-1}\omega \epsilon_4)$ and
$\nu_\alpha(\epsilon_5\eta)$ (which may be $<,=$ or $>$) belongs
to the list $h_1,\ldots,h_{p-1}$. Therefore, this relation is the same
for $\alpha_0$ and $\beta_0$. From this, it follows that the relation
between $\nu_\delta(Q_4)-\nu_{\delta 0}(Q_4)=\nu_\delta(Q'_4)$ and
$\nu_\delta(Q_5)-\nu_{\delta 0}(Q_5)=\nu_\delta(Q'_5)$ (which may be $<,=$ or
$>$) is
the same for $\delta=\alpha_0$ and for $\delta=\beta_0$. This
contradicts (\ref{eqn-compare1}), (\ref{eqn-compare2}).
$\Box$
\medskip

We now need the following geometric lemma.

\begin{lem} - \textbf{Tetrahedron Lemma} - Consider 4 points of $\R^3$,
$P,Q,R,S$, not lying in the same plane. Then they
define an affine basis of $\R^3$ and let $(u,v,w,t)$ be the
barycentric coordinates with respect to this basis (so
$u+v+w+t=1$).  Let $A$ be a point with coordinates
$(u_1,v_1,w_1,t_1)$ such that $u_1,v_1,w_1,t_1> 0$ and $u_1=v_1 \leq
w_1$ and $B$ be a point with coordinates $(u_2,v_2,w_2,t_2)$ such that
$u_2,v_2,w_2,t_2> 0$ and $v_2=w_2 \leq u_2$. Consider a finite set of
linear inequalities of the form $g_i \geq 0$, $i=1,\ldots,m$ such that
for all $i\in \{1,\ldots,m\}$, $g_i(A) \geq 0$ and $g_i(B) \geq
0$. Moreover suppose that, for each given $i$, only one of $u,v,w$
appears in $g_i$, which means that the plane $g_i=0$ passes at least
through two of the points $P,Q,R$.

Then there exists a point $D$ with coordinates
$(\lambda/3,\lambda/3,\lambda/3,1-\lambda)$, $0 \leq \lambda \leq 1$
such $g_i(D) \geq 0$ for all $i \in \{1,\ldots,m\}$.
\end{lem}
\begin{figure}[!h]
\begin{center}
\includegraphics[scale=0.9]{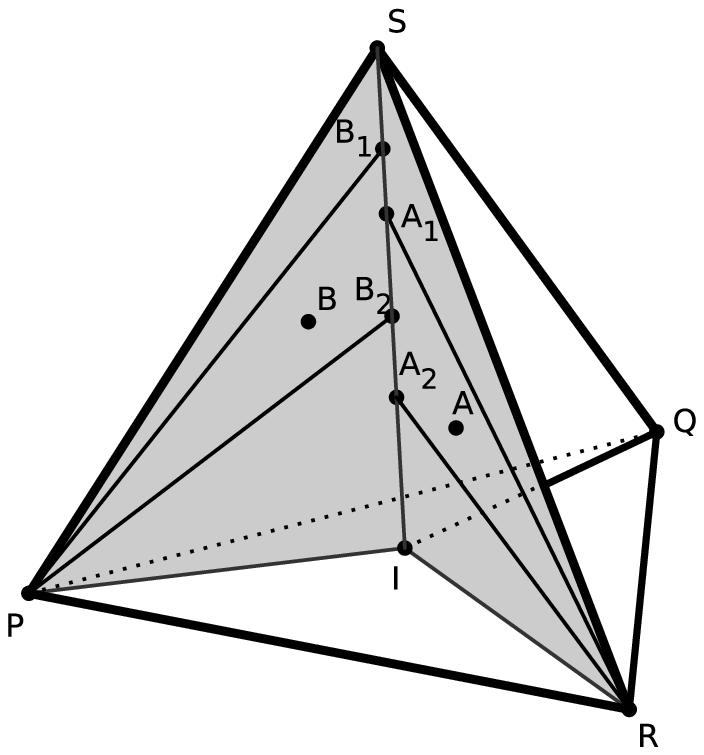}
\end{center}
\end{figure}
\noi Proof : Let $I$ be the point with coordinates
$(1/3,1/3,1/3,0)$. Without loss of generality, we may assume that, for
all $i$, the plane $\{g_i=0\}$ intersects the segment $IS$ in a point $M_i$
between $I$ and
$S$.

Let $i$ be such that $g_i(R)=0$. Then the intersection of $\{g_i=0\}$ with
the triangle $RIS$ is the segment $RM_i$.

Let $A_1$ and $A_2$ be the points of $IS$ defined by $$A_1S=
\bigcup_{\stackrel{RM_i \text{ lies above } A}{i \in \{1,\ldots,m\}}}M_iS$$ and
$$IA_2= \bigcup_{\stackrel{RM_i \text{ lies below } A}{i \in
\{1,\ldots,m\}}}IM_i$$ respectively.

Similarly, for all $i$ such that $g_i(P)=0$ we define $B_1,B_2 \in IS$
such that $$B_1S= \bigcup_{\stackrel{PM_i \text{ lies above } B}{i \in
\{1,\ldots,m\}}}M_iS $$ and $$IB_2= \bigcup_{\stackrel{PM_i \text{ lies below }
B}{i \in \{1,\ldots,m\}}}IM_i $$ respectively.

Let $B'$ be the point of intersection of $PB$ with $IS$.
We will prove that $A_1$ belongs to the interval $[SB']$ and $B'$ to $[SA_2]$.
And then put $D=B'$. This will complete the proof.

Let us show that $B' \in [SA_2]$. Now, if $A_2$ is defined by a plane of the
form $g_i=0$ with $g_i(R)=0$ and $g_i(P)=0$, then necessarily, $B'$ is closer
to $S$ than $A_2$.

So we may assume that $A_2$ is defined by a plane of the form
$u-kt=0$.

If $A_2$ is
closer to $S$ than $B'$, then $P$ and $B'$ are both in the same
half-space whose boundary is the plane $A_2QR$ with equation
$u-kt=0$. Then $(u-kt)(B) >0$. But, because,
$(u-kt)(A) < 0$ so also $(u-kt)(B) <0$, which is a contradiction. We prove the
same way that $A_1 \in [SB']$. So $B'$ lies between $A_1$ and $A_2$ as desired.
$\Box$
\medskip

Applying the Tetrahedron Lemma to all the relations of the form
(\ref{eqn-array12}) or (\ref{eqn-array22}) ensures the existence of a
point $D$ on $IS$ satisfying the same relations. Consider a point
$\delta \in \sper\ A$ such that $\nu_\delta$ has the same three
approximate roots $Q_4,Q_5,Q_6$ and satisfying moreover the fact that
the coordinates $(\nu(u_1),\nu(Q_4'),\nu(Q_5'),\nu(Q_6'))$ correspond to the
point $D$, then $\delta \in C$. So we are reduced to the limit case (see Remark
(\ref{rek_limit})) where
all the approximate roots are pairwise comparable for the couple
$(\alpha,\delta)$ and for the couple $(\beta,\delta)$. This defines two
connected sets $C_1$ and $C_2$, avoiding all $\{g_i=0\}$, and such that $C_1$
contains $\alpha$ and $\delta$ and $C_2$ contains $\beta$
and $\delta$. So letting $C=C_1 \cup C_2$ gives a connected set as
required.
\medskip

This settles the last remaining case ($(Q_4,Q_5,Q_6)$
pairwise incomparable) and with it Theorem \ref{th:ppal}.

\end{document}